\documentclass[preprint,review,12pt]{elsarticle}

\usepackage{amssymb}
\usepackage{amsmath}
\usepackage{amsthm}
\newtheorem{theorem}{Theorem}[section]

\numberwithin{equation}{section}

\newtheorem{lemma}[theorem]{Lemma}

\theoremstyle{definition}
\newtheorem{notation}[theorem]{Notation}
\newtheorem{remark}[theorem]{Remark}
\newtheorem{definition}[theorem]{Definition}
\newtheorem{normalization}[theorem]{Normalization}

\journal{Journal of Differential Equations}

\begin{document}

\begin{frontmatter}

\title{An Optimal Geometric Condition on Domains for Boundary Differentiability of Solutions of Elliptic Equations\tnoteref{t1}}

\tnotetext[t1]{Research supported by NSFC 11171266.}
\author[rvt]{Dongsheng Li}
\ead{lidsh@mail.xjtu.edu.cn}
\author[rvt]{Kai Zhang\corref{cor1}}
\ead{zkzkzk@stu.xjtu.edu.cn}

\cortext[cor1]{Corresponding author.}

\address[rvt]{School of Mathematics and Statistics, Xi'an Jiaotong University, Xi'an 710049, China}

\begin{abstract}
In this paper, a geometric condition on domains will be given which guarantees the boundary differentiability of solutions of elliptic equations, that is, the solutions are differentiable at any boundary point. We will show that this geometric condition is optimal.
\end{abstract}

\begin{keyword}
 Elliptic equations \sep Boundary differentiability \sep Optimal geometric condition
%% or \MSC[2008] code \sep code  (2000 is the default)

\end{keyword}

\end{frontmatter}

% \linenumbers

\section{Introduction}
In this paper, we will study the boundary differentiability of solutions of the following equations:
\begin{equation}\label{e1.1}
 \left\{\begin{aligned}
 -a^{ij}(x)\frac{\partial^2u(x)}{\partial x_i \partial x_j}=f(x)~~~~\mbox{in}~~\Omega; \\
 u=0~~~~\mbox{on}~~\partial\Omega,
\end{aligned}
\right.
\end{equation}
where $\Omega\subset R^n~(n>1)$ is a bounded domain; the matrix  $(a^{ij}(x))_{n\times n}\in C(\bar{\Omega})$ is symmetric and satisfies the uniformly elliptic condition with some constant $0<\lambda\leq1$, i.e., for any  $x\in \Omega$,
\begin{equation*}\label{e1.1.}
\lambda I_n\leq(a^{ij}(x))_{n\times n}\leq\frac{1}{\lambda}I_n
\end{equation*}
in the sense of nonnegative definiteness and $f\in C(\bar{\Omega})$. For convenience, solutions in this paper will always indicate viscosity solutions.

It is well known that the geometric properties of domains have significant influence on the boundary regularity of solutions. There are many remarkable results in this respect. First, Dirichlet problem asked what conditions on domains guarantee that the solutions of Laplace's equations are continuous at boundary. This problem was completely solved by Wiener \cite{Wi} in 1924 where the conception "Wiener Criterion" was introduced to describe such sufficient and necessary geometric properties of domains' boundaries. Besides, if $\Omega$ satisfies a uniform exterior cone condition, the solutions of (1.1) are H\"{o}lder continuous at the boundary (see Corollary 9.29\cite{GT}). Also, Trudinger obtained the Lipschitz continuity at the boundary under the hypothesis that $\Omega$ satisfies a uniform exterior sphere condition (see \cite{Tr} and \cite{Tr2}).

With regarding to the differentiability at boundary, Krylov got the $C^{1,\alpha}$ boundary regularity when $\partial \Omega$ belongs to $C^{1,\alpha}$ (see \cite{Kr} and \cite{Kr2}). Lieberman in \cite{Li} gave a more general result (Theorem 5.5) which contains above results in\cite{Kr2}. Furthermore, some results concerning Dini continuity can be found in \cite{Bu}, \cite{Il}, \cite{ME} and \cite{Sp}.

Almost all the previous results state that if the data are sufficiently smooth, then the solutions are $C^{1,\alpha}$ or $C^{1,Dini}$.  Li and Wang in \cite{LW} removed the smoothness assumptions on the boundary and obtained that the solutions of (1.1) are differentiable at boundary when $\Omega$ is convex. Later, the same authors got the boundary differentiability concerning the inhomogeneous boundary data condition (see \cite{LW2}). It is natural to ask whether convexity is the optimal geometric condition to guarantee the boundary differentiability. Actually, it is not and we will show that the boundary differentiability of solutions holds for a more general class of domains, which we call $\gamma-$convex domains. From Counterexample 4.1 and 4.3 in \cite{LW} and Theorem 1.11 in \cite{Sa}, we see that our result is optimal , i.e., the condition on the domain can not be weakened(cf. Remark \ref{r1.2}).

First, we define the concept $\gamma-$convexity.

\begin{definition}[\textbf{$\gamma$-convexity}]\label{d1.1}
Suppose that $ \Omega\subset R^n$ is a domain with continuous boundary. We call it $\gamma-$convex if the following holds:\\
There exist a constant $R_0(>0)$ and a function $\gamma:R^+\rightarrow R^+$, which is nondecreasing and satisfies
\begin{equation}\label{e1.3}
    \int_{0}^{R_0}\frac{\gamma(r)}{r^2}dr<\infty,
\end{equation}
such that for any $x_0\in\partial \Omega$ there exists a unit vector $\eta\in R^n$ such that
\begin{equation}\label{e1.2}
    \eta\cdot(x-x_0)\geq -\gamma(|x-x_0|),~\forall ~x\in\bar \Omega\cap B(x_0,R_0).
\end{equation}
\end{definition}

\begin{remark}\label{r1.1}
(i) Clearly, if $ \Omega$ is convex, then it is $\gamma-$convex with $\gamma\equiv 0$. Hence, results in \cite{LW} and \cite{LW2} are special cases of our results.\\
(ii) From (\ref{e1.3}), we see that $\frac{\gamma(r)}{r}\rightarrow 0$ as $r\rightarrow 0$. It follows that $\gamma(|x|)$ (as a function in $R^n$) is differentiable at 0 with $\gamma'(0)=0$. These properties will be used later.\\
(iii) If $x_0=0$ and $\eta=(0,0,...,0,1)$, then
\begin{equation*}\label{ee1.1}
    x_n\geq -\gamma(|x|),~\forall ~x\in\bar \Omega\cap B_{R_0}.
\end{equation*}
It follows that there exists $r_0>0$ (depending only on $\gamma$) such that
$$\frac{x_n}{|x'|}>-1 \mbox{ in } \bar\Omega\cap B_{r_0}.$$
Therefore,
\begin{equation}\label{ee1.2}
    x_n\geq-\gamma(2|x'|),~\forall ~x\in\bar \Omega\cap B_{r_0}.
\end{equation}
This means that there exists a differentiable hypersurface who touches $\partial\Omega$ at $0$ by below locally. This is the geometric explanation of $\gamma-$convexity.\qed

\end{remark}

%\begin{definition}\
%Let $\Omega\subset R^n$ be a domain. We say $\partial\Omega$ is continuous if for each point $x_0\in\partial\Omega$ there exist $r>0$ and a continuous function $\Gamma : R^{n-1}\rightarrow R$ such that
%\begin{equation*}
%    \Omega\cap B(x_0,r)=\{x|x_n>\Gamma(x_1,...,x_{n-1})\}\cap B(x_0,r),
%\end{equation*}
%upon relabeling and reorienting the coordinates axes if necessary.
%\end{definition}
Now we state our main result:

\begin{theorem}\label{t1.1}
Suppose that $\Omega$ is $\gamma-$convex and $u$ is the solution of (\ref{e1.1}). Then $u$ is differentiable at any boundary point. That is, for any $x_0\in \partial \Omega$, there exists a vector $a$ such that $u(x)=u(x_0)+a\cdot (x-x_0)+o(|x-x_0|),~\forall x\in\bar\Omega$.
\end{theorem}
\begin{remark}\label{r1.2}
Theorem \ref{t1.1} is optimal in the following two senses:

(1) $\gamma-$convexity can only guarantee the boundary differentiability and no more regularity can be expected. Li and Wang gave two counterexamples in \cite{LW} to show that the gradients of the solutions are not continuous.

(2) It is worth noting that to guarantee the boundary differentiability, the $\gamma-$convexity condition can not be weakened. Actually, Safonov proved the following result (see Theorem 1.11\cite{Sa}):

\emph{Let $x_0\in\partial \Omega$ and $u\in C^2(\Omega)\cap C(\bar\Omega)$ be a positive supersolution of}
\begin{equation*}
 -a^{ij}(x)\frac{\partial^2u(x)}{\partial x_i \partial x_j}=0~~~~\mbox{in}~~\Omega
\end{equation*}
\emph{such that $u=0$ on $\partial \Omega \cap B_{R_0}(x_0)$ for some $R_0>0$. }

\emph{If there exist a unit vector $\eta\in R^n$ and a function $\gamma:R^+\rightarrow R^+$, which is nondecreasing and satisfies}
\begin{equation*}\label{n14}
    \int_{0}^{R_0}\frac{\gamma(r)}{r^2}dr=\infty,
\end{equation*}
\emph{such that
$$\eta\cdot(x-x_0)\leq -\gamma(|x-x_0|),~~\forall ~x\in\bar {(\Omega)}^c\cap B(x_0,R_0),$$
then}
\begin{equation*}
\liminf_{t\rightarrow0^+}t^{-1}u(x_0+tl)=\infty~~for~all ~~l\cdot \eta>0.
\end{equation*}

Taking $f\equiv1$ in (\ref{e1.1}), then the classical solution $u$ of (\ref{e1.1}) satisfies the hypothesis above. Hence, $u$ can not be differentiable at $x_0$. Therefore, we see that the $\gamma-$convexity condition, thus (\ref{e1.3}) can not be released. From this point of view, our result is optimal.\qed
\end{remark}

This paper is organized as follows: The geometric properties of $\gamma-$convex domains will be studied in Section 2 where the main property is that blowing up at any boundary point of a $\gamma-$convex domain, we will obtain a cone. This is the same as convex domains. In Section 3, we will prove Theorem \ref{t1.1} by an iteration method, where the Harnack inequality and Aleksandrov-Bakelman-Pucci maximum principle are the main tools. We use the following notations in this paper, many of which are standard.

\begin{notation}
%%{\setlength{\baselineskip}{0.5\baselineskip}
\begin{enumerate}
  \item For any $x\in R^{n}$, we may write $x=(x',x_{n})$, where $x'\in R^{n-1}$ and $x_{n}\in R$.
  \item $\{e_i\}^{n}_{i=1}$: the standard basis of $R^n$.
  \item $R^n_+:=\{x\in R^n|x_n>0\}.$
  \item $ |x|:=\sqrt{\sum_{i=1}^{n} x_{i}^{2}},~\forall~ x\in R^{n}$.
  \item $A^c:$ the complement of $A$ in $R^n,~\forall A\subset R^n$.
  \item $\bar A: $ the closure of $A,~\forall~A\subset R^n$.
  \item $\mbox{dist}(A,B):=\inf\{|x-y||x\in A, y\in B\}$.
  \item $B(x_0,r):=\{x\in R^{n}|\ |x-x_0|<r\}$ \mbox{ and } $B_r:=B(0,r)$.
  \item $T_r\ :=\{x'\in R^{n-1}|\ |x'|<r\}$.
  \item $a^{+}:=\max\{0,a\}$ \mbox{ and } $a^{-}:=-\min\{0,a\}$.
  \item $Q[a,b]:=\{(x',x_n)\in R^{n}| \ |x'|<a, \ -b<x_n<b\}$.
  \item $Q_r:=Q[r,r]$.
  \item $\Omega[a,b]:=Q[a,b]\cap\Omega.$
  \item $\Omega_{r}:=\Omega[r,r]$.
\end{enumerate}
\end{notation}
\section{Geometric properties of $\gamma-$convex domains}

In this section, we study the blow-up sets at boundary points of $\gamma-$convex domains. It is well known that for a convex domain, the blow-up set at any boundary point is a cone. For a $\gamma-$convex domain, we will show the same conclusion.
\begin{definition}\label{d2.1}
Let $\Omega\subset R^n$ be a domain and $x_0\in \partial \Omega$. The blow-up set at $x_0$ is defined by
$$C_{x_0}=\left\{x|~\exists~t_x>0~\mbox{such that}~x_0+tx\in \Omega,~\forall ~0<t<t_x\right\}.$$
\end{definition}

\begin{definition}\label{d2.2}
We call $C\subset R^n$ a cone if it satisfies:\\
(i) If $x\in C$, then $t x\in C,~\forall~t >0$;\\
(ii) If $x,~y\in C$, then $x+y\in C$.\\
\end{definition}

The following lemma gives the relation between the boundary points of the blow-up set and that of the $\gamma-$convex domain.
\begin{lemma}\label{l2.1}
Suppose that $\Omega$ is $\gamma-$convex and $x_0\in\partial\Omega$. Then for any $\tilde x_0\in \partial C_{x_0}$ there exist a sequence $\{t_m\}$ monotone decreasing to 0, a sequence $\{z_m\}\subset\partial \Omega$, a unit vector sequence $\{\eta_m\}\subset R^n$ and a unit vector $\xi$ with $\eta_m\rightarrow\xi~\mbox{in}~R^n$ such that:

\noindent(i) $|x_0+t_m\tilde x_0-z_m|<t_m,~\forall~m\geq 1$;

\noindent(ii) $\eta_m\cdot(z-z_m)\geq -\gamma(|z-z_m|),~\forall~z\in \bar\Omega\cap B(z_m,R_0),~\forall~m\geq 1$;

\noindent(iii)  $\xi\cdot(\tilde x-\tilde x_0)\geq 0,~\forall~ x\in \overline C_{x_0}$,\\
where $\gamma$ and $R_0$ are the function and the constant in Definition \ref{d1.1}.
\end{lemma}

\noindent\textbf{Proof.}1. Without loss of generality, we assume that $x_0=0$.

2.  Suppose $\tilde x_0\not\in C_0$. Since $\tilde x_0\in\partial C_0$, there exists $\tilde y\in C_0$ such that $|\tilde x_0-\tilde y|<1$. By Definition \ref{d2.1}, there exists a sequence $\{t_m\}_{m=1}^{\infty}$ monotone decreasing to 0 such that $t_m\tilde x_0\in \Omega^c$ and there exists $t_{\tilde y}>0$ such that $t\tilde y\in \Omega,~\forall~0<t<t_{\tilde y}$. Without loss of generality, we assume that $t_1<t_{\tilde y}$. Thus, $t_m\tilde y\in \Omega$ for any $m\geq 1$. Then, for any $m\geq 1$, there exists $z_m\in\partial \Omega$ which lies in the line segment from $t_m\tilde x_0$ to $t_m\tilde y$. Hence,

\begin{equation*}\label{n3}
    |t_m\tilde x_0-z_m|\leq t_m|\tilde x_0-\tilde y|<t_m.
\end{equation*}
That is, (i) holds.

For any $m\geq 1$, since $ \Omega$ is $\gamma-$convex, there exists a unit vector $\eta_m\in R^n$ such that
\begin{equation*}\label{2}
\eta_m\cdot(z-z_m)\geq -\gamma(|z-z_m|),~ \forall ~z\in\bar \Omega\cap B(z_m,R_0).
\end{equation*}
That is, (ii) holds.

Since $\{\eta_m\}$ is bounded in $R^n$, there exists a subsequence of $\{\eta_m\}$, which we also denote by $\{\eta_m\}$, and a unit vector $\xi\in R^n$ such that
\begin{equation*}\label{3}
\eta_m\rightarrow\xi~\mbox{in}~R^n.
\end{equation*}

By Definition \ref{d2.1}, for any $\tilde x\in C_0$, there exists $t_{\tilde x}>0$ such that $t\tilde x\in \Omega,~\forall~0<t<t_{\tilde x}$. Without loss of generality, we assume that $t_1<t_{\tilde x}$. Thus, $t_m\tilde x\in \Omega$ for any $m\geq 1$. Hence,
\begin{eqnarray*}\label{n2}
    |t_m\tilde x-z_m|=|t_m\tilde x-t_m\tilde x_0+t_m\tilde x_0-z_m|\leq t_m(|\tilde x-\tilde x_0|+|\tilde x_0-\tilde y|).
\end{eqnarray*}
Therefore,
\begin{equation*}
    |t_m\tilde x-z_m|< R_0~\mbox{and}~|t_m\tilde x_0-z_m|< R_0,
\end{equation*}
for $m$ large enough.

Then, we have that
\begin{eqnarray*}
  \eta_m\cdot(t_m\tilde x-t_m\tilde x_0)&=&\eta_m\cdot(t_m\tilde x-z_m)+\eta_m\cdot(z_m-t_m\tilde x_0) \\
  &=&\eta_m\cdot(t_m\tilde x-z_m)+\lambda_m\eta_m\cdot(t_m\tilde y-z_m),
\end{eqnarray*}
where $\lambda_m = \frac{|z_m-t_m\tilde x_0|}{|t_m\tilde y-z_m|}$.
Hence,
\begin{eqnarray*}
   \eta_m\cdot(\tilde x-\tilde x_0)&=&\frac{1}{t_m}\big(\eta_m\cdot(t_m\tilde x-z_m)+\lambda_m\eta_m\cdot(t_m\tilde y-z_m)\big)\\
   &\geq&\frac{1}{t_m}\big(-\gamma(|t_m\tilde x-z_m|)-\gamma(|t_m\tilde x_0-z_m|)\big)\\
   &\geq&-\frac{1}{t_m}\Big(\gamma\big((|\tilde x-\tilde x_0|+|\tilde x_0-\tilde y|)t_m)\big)+\gamma(|\tilde x_0-\tilde y|t_m)\Big).\\
\end{eqnarray*}
Let $m\rightarrow\infty$, combining with Remark \ref{r1.1}(ii), we obtain that
$$\xi\cdot(\tilde x-\tilde x_0)\geq 0.$$
That is, (iv) holds.

3. Suppose $\tilde x_0\in C_0$. Then there exists $\tilde y\not\in C_0$ such that $|\tilde x_0-\tilde y|<1$.
Interchanging $\tilde x_0$ and $\tilde y$ in step 2, we can obtain $\{t_m\},~\{z_m\},~\xi~\mbox{and}~\{\eta_m\}$ such that (i), (ii) and (iii) hold similarly.\qed
~\\

The following theorem is the main result of this section which is a direct consequence of Lemma \ref{l2.1}.
\begin{theorem}\label{t2.1}
If $\Omega\subset R^n$ is $\gamma-$convex, then $C_{x_0}$  is a cone for any $x_0\in\partial \Omega$.
\end{theorem}
\noindent\textbf{Proof.} Without loss of generality, we assume that $x_0=0$. Let $\gamma$ and $R_0$ be the function and constant in Definition \ref{d1.1}. By Definition \ref{d2.1}, it is easy to see that
\begin{equation*}
\mbox{if}~\tilde x\in C_0, ~\mbox{then}~ t \tilde x\in C_0,~\forall~t>0.
\end{equation*}
Therefore, to show that $C_0$ is a cone, we only need to prove the convexity of $C_0$, which is equivalent to
\begin{equation}\label{e2.1}
\forall~ \tilde x_0\in\partial C_0,~\exists~\xi\in R^n~\mbox{such that}~\xi\cdot(\tilde x-\tilde x_0)\geq 0,~\forall~ \tilde x\in C_0.
\end{equation}
Let $\xi$ be given by Lemma \ref{l2.1}. By Lemma \ref{l2.1}(iii), (\ref{e2.1}) holds clearly.\qed

\begin{definition}\label{d2.3}
Let $\Omega$ be $\gamma-$convex. For any $x_0\in\partial \Omega$, if $\bar C_{x_0}\cap\bar{(-C_{x_0})}\neq R^{n-1}$, we call it a corner point. Otherwise, we call it a flat point.
\end{definition}

Before the end of this section, we prove the following two lemmas which give the geometric properties of the boundary with respect to corner points and flat points. These properties will be used to prove the boundary differentiability in the next section.
\begin{lemma}\label{t2.2}
Let $x_0\in\partial \Omega$ be a corner point and $\eta$ be the unit vector given in Definition \ref{d1.1}. Then
\begin{equation}\label{e2.2}
\exists~ y_0\in \Omega^c\mbox{ such that }\eta\cdot (y_0-x_0)>0~\mbox{and}~x_0+t(y_0-x_0)\in\Omega^c,~\forall~0<t<1.
\end{equation}
\end{lemma}

\noindent\textbf{Proof.} 1. Without loss of generality, we assume that $x_0=0$ and $\eta=(0,0,...,0,1)$. Thus $C_0\subsetneqq R^n_+$. Let $\gamma$ and $R_0$ be the function and constant in Definition \ref{d1.1}. Instead of proving (\ref{e2.2}), we only need to prove the following:

\begin{equation}\label{e2.3}
    \exists~ y\in R^n_+,~\exists~ t_y>0\mbox{ such that }ty\in\Omega^c,~\forall~0<t<t_y.
\end{equation}

2. Now, we choose the vector $y$ we need. Since $C_0\subsetneqq R^n_+$, there exists $\tilde x_0\in\partial C_0$ with $\tilde x_{0n}>0$ (we write $\tilde x_0=(\tilde x_{01},\tilde x_{02},...,\tilde x_{0n})$). Let $\{t_m\}$, $\{z_m\}$, $\eta_m$ and $\xi$ be given by Lemma \ref{l2.1} (with $x_0=0$). Let
\begin{equation*}
\xi^1=(\lambda_1\xi',\cos(\frac{1}{2}\arccos\xi_n)),
\end{equation*}
where $\lambda_1>0$ is chosen such that $|\xi^1|=1$. From $0\in \overline C_0$ and Lemma \ref{l2.1}(iii), we have $\xi\cdot\tilde x_0=0$. Thus, combining with $\tilde x_{0n}>0$ and $|\xi|=1$, we have $\xi_n<1$. Hence, $\xi_n<\xi^1_n<1$. Then, it is easy to see that $\xi^1\cdot\tilde x_0>0$. Without loss of generality, we may assume that $\xi^1\cdot\tilde x_0>2$ (otherwise, we may consider $\frac{3\tilde x_0}{\xi^1\cdot\tilde x_0}$). Combining with Lemma \ref{l2.1}(i), we have that
\begin{equation}\label{e2.5}
    \xi^1\cdot z_m>0,~\forall~m\geq 1.
\end{equation}

Since $\xi^1_n<1$, there exists $y\in R^n_+$ such that
\begin{equation}\label{e2.6}
    \xi^1\cdot y<0.
\end{equation}

3. Suppose (\ref{e2.3}) not hold. That is, there exists a sequence $\{\Tilde t_m\}$ monotone decreasing to 0 such that $\Tilde t_my\in\Omega$. From (\ref{e2.5}) and (\ref{e2.6}), we have that
\begin{equation*}
    \xi^1\cdot (\Tilde t_my-z_m)<0,~\forall~m\geq 1.
\end{equation*}
Hence,
$$\limsup_{m\rightarrow\infty} \frac{\xi^1\cdot (\Tilde t_my-z_m)}{|\Tilde t_my-z_m|}\leq 0.$$
On the other hand, since $|z_m|\rightarrow 0$, we may assume that the $n^{th}$ exponent of $(\Tilde t_my-z_m)$ is positive. Combining with $\xi^1_n>\xi_n$, we have that
\begin{equation*}
    \xi\cdot (\Tilde t_my-z_m)<0.
\end{equation*}
Hence,
$$\limsup_{m\rightarrow\infty} \frac{\xi\cdot (\Tilde t_my-z_m)}{|\Tilde t_my-z_m|}<0.$$

Finally,
$$\eta_m\cdot(\Tilde t_my-z_m)=(\eta_m-\xi)\cdot(\Tilde t_my-z_m)+\xi\cdot(\Tilde t_my-z_m).$$
Divide by $|\Tilde t_my-z_m|$ in above equation and let $m\rightarrow\infty$. Consequently, from Lemma \ref{l2.1}(ii) and (iii), the left side $\geq 0$ but the right side $<0$. We obtain a contradiction. Thus, (\ref{e2.3}) holds.\qed

\begin{lemma}\label{l2.3}
If $x_0\in\partial \Omega$ is a flat point, then $\partial \Omega$ is differentiable at $x_0$.
\end{lemma}

\noindent\textbf{Proof.} 1. Without loss of generality, we assume that $x_0=0$, $\eta=(0,...,0,1)$ and $\{x|x'=0,0<x_n\leq 1\}\subset \Omega$ where $\eta$ is the unit vector in Definition \ref{d1.1}. Hence, $\overline {C_0}\equiv\overline{ R^n_+}$. Let $\gamma$ and $R_0$ be the function and constant in Definition \ref{d1.1}.

2. Since $(0,...,0,1)\in \Omega$, there exists $\bar R>0$ such that
\begin{equation*}\label{7}
    \{x||x'|\leq \bar R,~x_n=1\}\subset \Omega.
\end{equation*}
Since $\frac{\gamma(r)}{r}\rightarrow 0$ as $r\rightarrow 0$, there exists $\widehat{R}>0$ such that
\begin{equation*}\label{9}
    0\leq \gamma(|x-y|)\leq \frac{1}{4n}|x-y| ~\mbox{ whenever }~ |x-y|<\widehat R.
\end{equation*}
Since $R_+^n$ is the blow-up set at 0, there exist $R_1,R_2,...,R_{n-1}>0$ such that
\begin{equation*}\label{8}
    \sup\{\frac{x_n}{t}|(te_i,x_n)\in\partial \Omega,~|t|<R_i\}<\frac{1}{2},~\forall~i=1,2,...,n-1,
\end{equation*}
where $\{e_i\}_{i=1}^{n-1}$ are the standard basis in $R^{n-1}$. Let
$$R_0=\min\{1,\bar R,\widehat R,R_0,R_1,R_2,...,R_{n-1}\}.$$

3. Suppose that $\partial\Omega$ is not differentiable at 0. That is,
$$\exists~\epsilon_0>0,~\forall~r>0,~\exists ~r'<r,~\exists~x\in \partial \Omega~\mbox{with}
~|x'|=r'~\mbox{such that}~\frac{x_n}{r'}\geq\epsilon_0.$$
We may assume that $\epsilon_0<1$ and that there exist $\tilde R<\frac{R_0}{4}$ and $y\in\partial \Omega$ such that $y_n=\epsilon_0\tilde R$ and $|y'|=\tilde R$. By a translation on the coordinate system, we assume that $y=(0,0,...,0,\epsilon_0\tilde R)$. From step 2, we have that
\begin{eqnarray*}
\lefteqn{S=\{x^1=(\frac{\epsilon_0\tilde R}{2},0,0,...,0,\frac{\epsilon_0\tilde R}{2}),x^2=(0,...,0,2\epsilon_0\tilde R),}\\
&&y^1_+=(\frac{\epsilon_0\tilde R}{2},0,0,...,0,\epsilon_0\tilde R),y^1_-=(-\frac{\epsilon_0\tilde R}{2},0,0,...,0,\epsilon_0\tilde R),...,\\
&&y^{n-1}_+=(0,0,0,...,\frac{\epsilon_0\tilde R}{2},\epsilon_0\tilde R),y^{n-1}_-=(0,0,0,...,-\frac{\epsilon_0\tilde R}{2},\epsilon_0\tilde R)\}\subset \Omega.
\end{eqnarray*}

Since $y\in\partial \Omega$, there exists a unit vector $\tilde \eta=(\tilde \eta_1,\tilde \eta_2,...,\tilde \eta_n)$ such that
\begin{equation*}\label{10}
    \tilde \eta\cdot(x-y)\geq -\gamma(|x-y|)\geq -\frac{1}{4n}|x-y|,~~\forall~ x\in \bar\Omega\cap B(y,R_0).
\end{equation*}
Combining with $y^i_+\in \Omega\cap B(y,R_0)$, we have that
\begin{equation*}
 \tilde \eta_i\cdot\frac{1}{2}\epsilon_0\tilde R=\tilde \eta\cdot(y^i_+-y)\geq -\gamma (|y^i_+-y|)\geq -\frac{1}{4n}\cdot\frac{1}{2}\epsilon_0 \tilde R,~\forall~i=1,2,...,n-1.
\end{equation*}
Hence,
\begin{equation*}\label{n8}
    \tilde \eta_i\geq -\frac{1}{4n},~\forall~i=1,2,...,n-1.
\end{equation*}
Similarly, since $y^i_-\in \Omega\cap B(y,R_0)$, we have that
\begin{equation*}
    \tilde \eta_i\leq \frac{1}{4n},~\forall~i=1,2,...,n-1.
\end{equation*}
Since $x^1\in \Omega\cap B(y,R_0)$, we have that
\begin{equation*}
    \tilde \eta_n \leq \tilde \eta_1+\frac{\sqrt{2}}{8n}\leq \frac{1}{2n}.
\end{equation*}
Since $x^2\in \Omega\cap B(y,R_0)$, we have that
\begin{equation*}
    \tilde \eta_n \geq -\frac{1}{4n}.
\end{equation*}
Hence, we obtain a contradiction with $|\tilde \eta|=1$. Therefore, $\partial\Omega$ is differentiable at 0.\qed

\section{Differentiability at the boundary}

In this section we will prove Theorem \ref{t1.1}. Our proof will be divided into two parts according to the two kinds of boundary points:
corner points and flat points (cf. Definition \ref{d2.3}). In order to prove the theorem concisely, we make the following normalizations
without loss of generality.

\begin{normalization}\label{nor1}

(1) We assume that $0\in\partial\Omega$ and we only prove that $u$ is differentiable at $0$.

(2) By the linearity of the equation, we assume that $u,f\geq 0~~\mbox{in}~\Omega$ and $||u||_{L^{\infty}(\Omega_1)},||f||_{L^{n}(\Omega_1)}\leq 1$.

(3) We assume that $\eta=(0,...,0,1)$ and $R_0=2$ in Definition \ref{d1.1}.

(4) We may write $\gamma(r)= r\sigma(r)$ and assume that
$$\int_0^1\frac{\sigma(r)}{r}dr\leq 1.$$

(5) Since $\frac{\gamma(r)}{r}\rightarrow 0$ as $r\rightarrow 0$, there exists $r_1>0$ such that $\frac{\gamma(r)}{r}<C~(\forall~r\leq r_1)$ where $C$ depending only on $n$ and $\lambda$ is small enough (in fact, we may take $C=\frac{1}{32(n-1)(1+\frac{2\sqrt{n-1}}{\lambda})^2}$). We assume that $r_1=1$. Hence, we have that
$$\gamma(t)<tC,~\forall~ 0<t\leq 1.$$

(6) By (\ref{ee1.2}) in Remark \ref{r1.1}, we assume that
\begin{equation*}
    x_n\geq-\gamma(2|x'|),~\forall ~x\in\bar \Omega\cap B_{2}.
\end{equation*}
Combining with (5), we have that
\begin{equation}\label{ecee3.1}
    x_n\geq-\gamma(2|x'|),~\forall ~(x',x_n)\in\bar \Omega, \mbox{ where } x'\in T_{1}.
\end{equation}
\end{normalization}

In the following lemma, we construct two barrier functions which will be used later repeatedly.

\begin{lemma}\label{l3.1}
Suppose that $M\geq \sqrt{n-1}(1+\frac{4\sqrt{n-1}}{\lambda})$, $\delta >0$, $M\delta<\frac{1}{2}$ and $\gamma(2M\delta)<\delta$. Then there exist two twice differentiable functions $\Psi_{M\delta,\delta}$ and $\widetilde{\Psi}_{M\delta,\delta}$ which satisfy

\begin{eqnarray}\label{11}
 \left\{\begin{array}{l}
(1)~\Psi_{M\delta,\delta}\geq1~~~~~~if ~x_{n}=\delta~~\mbox{and}~~|x'|\leq M\delta,\\[3pt]
(2)~\Psi_{M\delta,\delta}\geq0~~~~~~\mbox{on}~~ \overline{\Omega[M\delta,\delta]},\\[3pt]
(3)~\Psi_{M\delta,\delta}\geq1~~~~~~if ~x\in\overline{\Omega[M\delta,\delta]}~~\mbox{and}~~|x'|=M\delta,\\[3pt]
(4)\displaystyle~\Psi_{M\delta,\delta}\leq\frac{2(x_n+\gamma(2M\delta))}{\delta}~~~~~~\mbox{on}~~\overline{
\Omega}_{\delta},\\[3pt]
(5)\displaystyle~-a^{ij}(x)\frac{\partial^2\Psi_{M\delta,\delta}(x)}{\partial x_i \partial x_j}\geq 0~~~~~~\mbox{in}~~\Omega[M\delta,\delta],\\[3pt]
\end{array}
\right.
\end{eqnarray}
and
\begin{eqnarray}\label{12}
 \left\{\begin{array}{l}
(1)~\widetilde{\Psi}_{M\delta,\delta}\leq1~~~~~~if~ x_{n}=\delta~~\mbox{and}~~|x'|\leq M\delta,\\[3pt]
(2)\displaystyle~\widetilde{\Psi}_{M\delta,\delta}\leq\frac{x_{n}+\gamma(2M\delta)}{\delta}~~~~~~\mbox{on}~~ \overline{\Omega[M\delta,\delta]},\\[3pt]
(3)~\widetilde{\Psi}_{M\delta,\delta}\leq0~~~~~~if ~x\in\overline{\Omega[M\delta,\delta]}~~\mbox{and}~~|x'|=M\delta,\\[3pt]
(4)\displaystyle~\widetilde{\Psi}_{M\delta,\delta}\geq\frac{x_n+\gamma(2M\delta)}{4\delta}~~~~~~\mbox{on}~~\overline{
\Omega}_{\delta},\\[3pt]
(5)\displaystyle~-a^{ij}(x)\frac{\partial^2\widetilde{\Psi}_{M\delta,\delta}(x)}{\partial x_i \partial x_j}\leq0~~~~~~\mbox{in}~~\Omega[M\delta,\delta],\\[3pt]
\end{array}
\right.
\end{eqnarray}
respectively.
\end{lemma}
\noindent\textbf{Proof.} Choose $\epsilon~>0$ small enough such that
\begin{equation}\label{e3.1}
4-(1+\epsilon)(2+\epsilon)(M-1)^{\epsilon}\geq0.
\end{equation}
Define $\Psi_{M\delta,\delta}$ and $\widetilde{\Psi}_{M\delta,\delta}$ as follows:
\begin{equation*}
\Psi_{M\delta,\delta}(x)=\frac{2(x_n+\gamma(2M\delta))}{\delta}-\left(\frac{x_n+\gamma(2M\delta)}{2\delta}\right)^{2}
+\frac{\lambda^2}{8(n-1)}\sum_{i=1}^{n-1}\left(\left(\frac{|x_i|}{\delta}-1\right)^{+}\right)^{2+\epsilon}
\end{equation*}
and
\begin{equation*}
\widetilde{\Psi}_{M\delta,\delta}(x)=\frac{1}{2} \left(\frac{x_n+\gamma(2M\delta)}{2\delta}+\left (\frac{x_n+\gamma(2M\delta)}{2\delta}\right)^{2}\right) -\frac{\lambda^2}{16(n-1)}\sum_{i=1}^{n-1}\left(\left(\frac{|x_i|}{\delta}-1\right)^{+}\right)^{2+\epsilon}.
\end{equation*}

Since the proofs of (\ref{11}) and (\ref{12}) are similar, we only prove (\ref{11}). By the definition of $\Psi_{M\delta,\delta}$, (\ref{11}(1)) and (\ref{11}(4)) hold clearly.

From (\ref{ecee3.1}), we have that
\begin{eqnarray*}
  x_n\geq -\gamma(2|x'|)\geq -\gamma(2M\delta)~~~~\mbox{on}~~\overline{\Omega[M\delta,\delta]}.
\end{eqnarray*}
Hence, (\ref{11}(2)) holds.

If $x\in\overline{\Omega[M\delta,\delta]}$, we have that
\begin{equation}\label{e3.2}
    0\leq \frac{x_n+\gamma(2M\delta)}{2\delta}\leq 1.
\end{equation}
If $|x'|=M\delta$, there exists an index $i$ such that
$$|x_i|\geq \frac{M\delta}{\sqrt{n-1}}.$$
Thus, from $M\geq \sqrt{n-1}(1+\frac{4\sqrt{n-1}}{\lambda})$ and (\ref{e3.2}), we have that
$$\Psi_{M\delta,\delta}\geq \frac{\lambda^2}{8(n-1)}\left(\frac{M}{\sqrt{n-1}}-1\right)^{2+\epsilon}\geq 2.$$
That is, (\ref{11}(3)) holds.

Finally, from (\ref{e3.1}), we have that
\begin{eqnarray*}
-a^{ij}(x)\frac{\partial^2\Psi_{M\delta,\delta}(x)}{\partial x_i \partial x_j}\geq \frac{a^{nn}}{2\delta^2}-\frac{(2+\epsilon)(1+\epsilon)\lambda^2}{8(n-1)\delta^2}\sum_{i=1}^{n-1}a^{ii}(M-1)^\epsilon
\geq \frac{\lambda}{2\delta^2}-\frac{4\lambda}{8\delta^2}
= 0.
\end{eqnarray*}
That is, (\ref{11}(5)) holds.\qed
~\\

\subsection{Differentiability at corner points}

In this subsection, we prove Theorem \ref{t1.1} with respect to corner points. By Lemma \ref{t2.2}, there exist $h_0>0$ and $r_0>0$ such that
\begin{equation*}
    th_0e_n+T_{tr_0}\cap\partial\Omega\neq\emptyset,~\forall~0<t<1.
\end{equation*}
Without loss of generality, we assume that $h_0=1$. Then combining with Normalization \ref{nor1}, Theorem \ref{t1.1} is a direct consequence of the following:
\begin{theorem}\label{ectt3.1}
There exist positive constants $C$, $\alpha$ and $\Lambda$ depending only on $n$, $\lambda$, $r_0$ and $\sigma$ such that
\begin{equation}
u(x)\leq Cr\bigg\{r^{\alpha}\bigg(1+\int_{r}^{1}\frac{\|f\|_{L^{n}(\Omega_{t})}+\sigma(t)}{t^{1+\alpha}}dt\bigg)+\|f\|_{L^{n}(\Omega_{\Lambda r})}+\sigma(\Lambda r)\bigg\},
\end{equation}
for any $x\in\Omega_r$ and $r\leq \frac{1}{\Lambda}$.
\end{theorem}

In this subsection, unless otherwise stated, $C$, $\tilde C$, $C_1$, $C_2$, etc., denote constants depending only on $n$, $\lambda$, $r_0$ and $\sigma$.

Now, we use an iteration method to prove Theorem \ref{ectt3.1} where the Aleksandrov-Bakelman-Pucci maximum principle and the Harnack inequality are the main tools.
\begin{lemma}[\textbf{Key iteration}]\label{ecl3.2}
There exist positive constants $\delta(<1)$, $\mu(<1)$, $M$, $A_{1}$ and $A
_{2}$ depending only on $n$, $\lambda$ and $r_0$ such that if
\begin{equation}\label{ece3.3}
u(x)\leq Kx_n+B~~~~in~~\Omega_1
\end{equation}
for some nonnegative constants $K$ and $B$, then
\begin{eqnarray}\label{ec14}
\begin{split}
u(x)\leq \mu(K+MB)x_n+A_{1}\|f\|_{L^{n}(\Omega_{1})}+A_{2}(K+B)\gamma(1)~~~~in~~\Omega_{\delta}.
\end{split}
\end{eqnarray}
\end{lemma}

\noindent\textbf{Proof.} 1. \emph{Claim:} There exist positive constants $M$, $\delta_1$ and $C_1$ depending only on $n$ and $\lambda$
such that
\begin{eqnarray}\label{ec17}
\begin{split}
u(x)\leq (K+MB)x_n+C_1\|f\|_{L^{n}(\Omega_{1})}+(K+MB)\gamma(1)~~~~\mbox{in}~~\Omega_{\delta_1}.
\end{split}
\end{eqnarray}
\proof Let $M_1=2\sqrt{n-1}\left(1+\frac{2\sqrt{n-1}}{\lambda}\right)$, $M=4M_1$, $\delta_1=\frac{1}{2M_1}$ and $\Psi=\Psi_{\frac{1}{2},\delta_1}$ (see Lemma \ref{l3.1} and recall $\gamma(1)\leq\frac{\delta_1}{4}$). Let $v(x)=u(x)-Kx_n-B\Psi(x)$. We claim that
\begin{equation*}
    v(x)\leq K\gamma(1)~~~~~~~\mbox{on}~~\partial\Omega[\frac{1}{2},\delta_1].
\end{equation*}
In fact, we separate $\partial\Omega[\frac{1}{2},\delta_1]$ into two parts: $\partial Q[\frac{1}{2},\delta_1]\cap\bar\Omega$ and $\partial\Omega\cap \overline{ Q[\frac{1}{2},\delta_1]}$. On the first part, since $\Psi\geq 1$, we have $v\leq 0$. On the second part, since $\Psi\geq 0$, $u=0$ and $x_n\geq -\gamma(1)$, we have $v\leq K\gamma(1)$. Thus, we have proved the claim.
Then it is easy to verify that
\begin{eqnarray*}
 \left\{\begin{array}{l}
 \displaystyle-a^{ij}(x)\frac{\partial^2v(x)}{\partial x_i \partial x_j}\leq f(x)~~~~\mbox{in}~~\Omega[\frac{1}{2},\delta_1];\\
v\leq K\gamma(1)~~~~\mbox{on}~~\partial\Omega[\frac{1}{2},\delta_1].\\
\end{array}
\right.
\end{eqnarray*}

According to the Aleksandrov-Bakelman-Pucci maximum principle\cite{GT}, we have that
\begin{equation*}
v(x)\leq K\gamma(1)+C_1||f||_{L^{n}(\Omega_1)}~~~~\mbox{in}~~\Omega[\frac{1}{2},\delta_1].
\end{equation*}
From (\ref{11}(4)) and the definition of $v$, we have that
\begin{eqnarray*}\label{ec18}
u(x)&=&Kx_n+B\Psi(x)+v(x)\\
&\leq& (K+\frac{2B}{\delta_1})x_n+C_1||f||_{L^{n}(\Omega_1)}+(K+\frac{2B}{\delta_1})\gamma(1)~~~~\mbox{in}~~\Omega_{\delta_1}.
\end{eqnarray*}
Since $M=\frac{2}{\delta_1}$, (\ref{ec17}) holds.

2. Let $M_1=\max\{M_1,r_0\}$, $\delta=\frac{\delta_1}{2M_1}$ and $\Gamma=(\delta e_n+T_{M_1\delta})\cap\Omega$ (recall $\gamma(1)<\frac{1}{4}\delta$). Then $\bar\Gamma\cap\partial\Omega\neq\emptyset.$ Let
\begin{equation*}
    \Gamma_0=\{x\in \Gamma | \mbox{dist}(x,\partial\Omega)\leq d\}\mbox{ and }\Gamma_1=\overline{\Gamma\backslash\Gamma_0},
\end{equation*}
where $d$ is a constant depending only on $n$ and $\lambda$ and is determined by the following way. Since $\Omega$ is $\gamma-$convex, $\Omega$ satisfies the uniform exterior cone condition. Therefore, $u$ is H\"{o}lder continuous on $\bar\Omega$. Then by Corollary 9.28 in \cite{GT}, we have that
\begin{equation}\label{ec3.1}
   |u(x)-u(y)|\leq \tilde C\bigg\{\Big(\frac{d}{\delta_1}\Big)^{\mu_0}(||u||_{L^{\infty}(\Omega_{\delta_1})}+||f||_{L^{n}(\Omega_{\delta_1})})\bigg\},
\end{equation}
for any $x\in\Gamma_0$ and $y\in\bar\Gamma_0\cap\partial\Omega$, where $\tilde C$ and $\mu_0$ are positive constants depending only on $n$ and $\lambda$. Let $d$ be small enough such that $\tilde C(\frac{d}{\delta_1})^{\mu_0}\leq \frac{1}{2}.$ Then combining with (\ref{ec17}), we have that
\begin{equation*}
    u(x)\leq \frac{1}{2}(K+MB)\delta+(C_1+1)\|f\|_{L^{n}(\Omega_{1})}+\frac{1}{2}(K+MB)\gamma(1),~\forall~x\in\Gamma_0.
\end{equation*}

Let
$$v(x)=(K+MB)x_n+C_1\|f\|_{L^{n}(\Omega_{1})}+(K+MB)\gamma(1)-u(x).$$
Then
\begin{eqnarray*}
  \inf_{\Gamma_0} v(x)&\geq& \frac{1}{2}(K+MB)\delta-\|f\|_{L^{n}(\Omega_{1})}+\frac{1}{2}(K+MB)\gamma(1).\\
\end{eqnarray*}
First, we assume that $\Gamma_1\neq\emptyset$. Thus, $\Gamma_0\cap\Gamma_1\neq\emptyset$ and thereby
\begin{equation*}
    \inf_{\Gamma_0}v(x)\leq \sup_{\Gamma_1}v(x).
\end{equation*}
Since $v\geq0~~\mbox{in}~~ \Omega_{\delta_1}$, we apply the Harnack inequality to $v$ on $\Gamma_1$ and have that
\begin{equation*}
    \sup_{\Gamma_1}v(x)\leq C_2\left(\inf_{\Gamma_1}v(x)+\|f\|_{L^{n}(\Omega_{1})}\right).
\end{equation*}
Hence,
\begin{equation*}\label{ec29}
  \inf_{\Gamma}v(x)\geq\left\{\frac{1}{2C_2}\Big((K+MB)\delta+(K+MB)\gamma(1)\Big)-\Big(\frac{1}{C_2}+1\Big)\|f\|_{L^{n}(\Omega_{1})}\right\}^{+}:=a.
\end{equation*}
Clearly, if $\Gamma_1=\emptyset$, we also have $\inf_{\Gamma}v(x)\geq a$.

Let $\tilde{\Psi}=\tilde\Psi_{M_1\delta,\delta}~\mbox{and}~w=a\tilde{\Psi}-v$. We claim that
\begin{equation*}\label{ec30}\
    w\leq 0~~~~~~\mbox{on}~~\partial\Omega[M_1\delta,\delta].
\end{equation*}
In fact, we separate $\partial\Omega[M_1\delta,\delta]$ into three parts:
$\{x\in\partial Q[M_1\delta,\delta]|x_n=\delta\}\cap\bar\Omega$, $\{x\in\partial Q[M_1\delta,\delta]|x_n\neq\delta\}\cap\bar\Omega$ and $\partial\Omega\cap \overline{Q[M_1\delta,\delta]}.$ On the first part, since $\tilde\Psi\leq 1$ and $v\geq a$, we have $w\leq 0$. On the second part, since $\tilde\Psi\leq 0$ and $v\geq 0$, we have $w\leq 0$. On the last part, since $\tilde{\Psi}(x)\leq\frac{x_n+\gamma(2M_1\delta)}{\delta}$, we have that
\begin{eqnarray*}
  a\tilde{\Psi}(x)-v(x)&\leq& \frac{1}{2C_2}\Big((K+MB)\delta+
    (K+MB)\gamma(1)\Big)\frac{x_n+\gamma(2M_1\delta)}{\delta}-v(x)\\
 &\leq& \frac{1}{2C_2}(K+MB)\Big(x_n+\gamma(2M_1\delta)\Big)+\frac{1}{2C_2}(K+MB)\gamma(1)
 \frac{x_n+\gamma(2M_1\delta)}{\delta}\\
 &&-\bigg((K+MB)\Big(x_n+\gamma(1)\Big)+C_1\|f\|_{L^{n}(\Omega_{1})}-u(x)\bigg)\\
 &\leq& u(x)=0
  ~~~~\mbox{on}~~ \partial\Omega\cap \overline{Q[M_1\delta,\delta]}.
\end{eqnarray*}
Thus, we have proved the claim. Then it is easy to verify that
\begin{eqnarray*}
 \left\{\begin{array}{l}
\displaystyle-a^{ij}(x)\frac{\partial^2w(x)}{\partial x_i \partial x_j}\leq f(x)~~~~\mbox{in}~~\Omega[M_1\delta,\delta];\\
w\leq 0~~~~\mbox{on}~~ \partial\Omega[M_1\delta,\delta].\\
\end{array}
\right.
\end{eqnarray*}

According to the Aleksandrov-Bakelman-Pucci maximum principle, we have that

\begin{equation*}
w\leq C_1\|f\|_{L^{n}(\Omega_{1})}~~~~\mbox{in}~~ \Omega[M_1\delta,\delta].
\end{equation*}
From (\ref{12}(4)) and the definition of $v$, we have that
\begin{eqnarray*}
  u(x)&=&(K+MB)x_n+C_1||f||_{L^{n}(\Omega_1)}+(K+MB)\gamma(1)-v(x) \\
  &=&(K+MB)x_n+C_1||f||_{L^{n}(\Omega_1)}+(K+MB)\gamma(1)-a\tilde\Psi(x)+w(x)\\
&\leq& (K+MB)x_n+C_1||f||_{L^{n}(\Omega_1)}+(K+MB)\gamma(1)-\\
&&\frac{a}{4\delta}\Big(x_n+\gamma(2M_1\delta)\Big)
+C_1\|f\|_{L^{n}(\Omega_{1})}\\
&\leq& (K+MB)x_n+2C_1||f||_{L^{n}(\Omega_1)}+(K+MB)\gamma(1)-\\
&&\frac{x_n+\gamma(2M_1\delta)}{4\delta}\left(\frac{K+MB}{2C_2}\delta-2||f||_{L^{n}(\Omega_1)}\right)\\
&=& (1-\frac{1}{8C_2})(K+MB)x_n+2(C_1+1)||f||_{L^{n}(\Omega_1)}+\\
&&(K+MB)\gamma(1)-\frac{\gamma(2M_1\delta)}{8C_2}(K+MB)\\
&\leq& (1-\frac{1}{8C_2})(K+MB)x_n+2(C_1+1)||f||_{L^{n}(\Omega_1)}+\\
&&(M+1)(K+B)\gamma(1)~~~~\mbox{in}~~\Omega_\delta.
\end{eqnarray*}
Let $\mu=1-\frac{1}{8C_2}$, $A_1=2(C_1+1)$, and $A_2=2M+1$. Then (\ref{ec14}) holds.\qed

\begin{lemma}\label{ecl3.3}
For $m=0,1,2,...$, let
\begin{equation*}
K_{m+1}=\mu(K_m+M\frac{B_m}{\delta^m})\mbox{~and~} B_{m+1}=A_1\delta^m||f||_{L^{n}(\Omega_{\delta^m})}+A_2(K_m+\frac{B_m}{\delta^m})\delta^m\sigma(\delta^m),
\end{equation*}
where $K_0=0$ and $B_0=1$. Then
\begin{equation}\label{ec8}
    u(x)\leq K_mx_n+B_m~~~~in~~ \Omega_{\delta^m},
\end{equation}
for $m=0,1,2,...$.
\end{lemma}

\noindent\textbf{Proof.} 1. \emph{Claim:} Let $\tilde\Omega=\frac{1}{t}\Omega~(t>0)$. Then $\tilde\Omega$ is also $\gamma-$convex with the function $\tilde\gamma(r)=\frac{1}{t}\gamma(tr)$ and the constant $\tilde R=\frac{2}{t}$ (see Definition \ref{d1.1}).

\proof Since $\Omega$ is $\gamma-$convex, we have that for any $x_0\in\partial \Omega$, there exists a unit vector $\eta\in R^n$ such that
$$\eta\cdot(x-x_0)\geq -\gamma(|x-x_0|)~~~~\forall ~x\in\bar \Omega\cap B(x_0,2).$$
Let $\tilde\eta=\eta$, we have that
\begin{eqnarray*}
\tilde\eta\cdot(\frac{1}{t}x-\frac{1}{t}x_0)
=\frac{1}{t}\eta\cdot(x-x_0)
\geq -\frac{1}{t}\gamma(|x-x_0|)
=-\frac{1}{t}\gamma(|t\frac{1}{t}x-t\frac{1}{t}x_0|)
=-\tilde\gamma(|\frac{1}{t}x-\frac{1}{t}x_0|).
\end{eqnarray*}
Hence $\tilde\Omega$ is $\gamma-$convex with the function $\tilde\gamma$ and the constant $\tilde R$.

2. We use induction method to prove the lemma. From Lemma \ref{ecl3.2} and $0\leq u \leq 1$ we see that (\ref{ec8}) holds for $m=0$. Suppose that (\ref{ec8}) holds for $m=l-1$.

Let $y=\frac{x}{\delta^l},~\tilde u(y)=\frac{u(\delta^ly)}{\delta^{2l}},~\tilde a^{ij}(y)=a^{ij}(\delta^ly),~\tilde f(y)=f(\delta^ly)~\mbox{and}~\tilde\Omega=\frac{1}{\delta^l}\Omega$. It is clear that $(a^{ij}(x))_{n\times n}$ satisfies the uniformly elliptic condition with $\lambda$ and $0\in\partial\tilde\Omega$. Let $\tilde\nu$ denotes the boundary of $\tilde\Omega$ near 0. Then $\tilde\nu$ satisfies that

\begin{equation*}
    \tilde\nu(y')=\frac{1}{\delta^l}\nu(\delta^ly')~~~~\mbox{where}~ y'\in T_1.
\end{equation*}
By Claim 1, $\tilde\gamma(1)=\frac{1}{\delta^l}\gamma(\delta^l)\leq \gamma(1)$. On the other hand, it is easy to verify that

\begin{eqnarray*}
 \left\{\begin{array}{l}
 \displaystyle -\tilde a^{ij}(y)\frac{\partial^2\tilde u(y)}{\partial y_i \partial y_j}=\tilde f(y)~~~~\mbox{in}~~\tilde\Omega; \\[10pt]
 \tilde u=0~~~~\mbox{on}~~\partial\tilde\Omega.
\end{array}
\right.
\end{eqnarray*}
From the induction assumptions and the definition of $\tilde u$, we have that

\begin{equation*}\label{ec35}
  \tilde u(y)\leq \frac{K_l}{\delta^l}y_n+\frac{B_l}{\delta^{2l}}~~~~~\mbox{in}~~\tilde\Omega_1.
\end{equation*}
Then, by Lemma \ref{ecl3.2}, we have that

\begin{eqnarray*}
    &&\tilde u(y)\leq \mu (\frac{K_l}{\delta^l}+M\frac{B_l}{\delta^{2l}})y_n+A_1||\tilde f||_{L^{n}(\tilde\Omega_1)}+
    A_2\Big(\frac{K_l}{\delta^l}+\frac{B_l}{\delta^{2l}}\Big)\tilde\gamma(1)~~~~\mbox{in}~\tilde\Omega_\delta,
\end{eqnarray*}
By variables changing, combining with $\tilde\gamma(1)=\frac{1}{\delta^l}\gamma(\delta^l)= \sigma(\delta^l)$, we have that (\ref{ec8}) hold for $m=l+1$.\qed
~\\

Before proving Theorem \ref{ectt3.1}, we first give the following two important facts which will be used both in
the proof of corner points and in the proof of flat points:
\begin{equation}\label{e3.11}
   \sum_{i=0}^{\infty}f_{\delta^{i}}<\infty~\mbox{and}~\sum_{i=0}^{\infty}\sigma_{\delta^{i}}<\infty.
\end{equation}
In fact
\begin{eqnarray*}
  \sum_{i=0}^{m}f_{\delta^{i}}=1+\frac{\delta}{1-\delta}\sum_{i=1}^{m}(\delta^{i-1}-\delta^{i})\frac{f_{\delta^{i}}}{\delta^i}
  \leq 1+\frac{\delta}{1-\delta}\int_{\delta^{m}}^{1}\frac{f_r}{r}dr.
\end{eqnarray*}
Since $f\in C(\bar\Omega)$, we obtain that $\sum_{i=0}^{\infty}f_{\delta^{i}}<\infty$. Similarly, $\sum_{i=0}^{\infty}\sigma_{\delta^{i}}<\infty$.
~\\

\noindent\textbf{Proof of Theorem \ref{ectt3.1}.} Let $\{B_m\}_{m=0}^{\infty}$ and
$\{K_m\}_{m=0}^{\infty}$ be defined as in Lemma \ref{ecl3.3}. For simplicity, we denote $||f||_{L^{n}(\Omega_r)}$ by $f_r$ and
$\sigma(r)$ by $\sigma_r$. For $m=2,3,...$, we have that

\begin{eqnarray*}
  K_m+\frac{B_m}{\delta^m} &=& \mu (K_{m-1}+M\frac{B_{m-1}}{\delta^{m-1}})+ \frac{B_{m}}{\delta^{m}}\\
   &=& M\Big(\mu^m+\sum_{i=1}^{m}\mu^i\frac{B_{m-i}}{\delta^{m-i}}\Big)+\frac{B_{m}}{\delta^{m}}\\
   &\leq& M\Big(1+\sum_{i=0}^{m}\frac{B_i}{\delta^i}\Big)\\
   &\leq& 2M\Big(1+\frac{1}{\delta}\sum_{i=1}^{m}\big(A_1f_{\delta^{i-1}}+A_2(K_{i-1}+\frac{B_{i-1}}{\delta^{i-1}})\sigma_{\delta^{i-1}}\big)\Big).
\end{eqnarray*}
From (\ref{e3.11}), it is easy to verify that
\begin{equation*}
    K_m+\frac{B_m}{\delta^m}\leq C_4.
\end{equation*}
Thus, by a simple calculation, we see that

\begin{equation*}
K_{m}\leq M\mu^{m}\bigg(1+\frac{A_1}{\delta}\sum^{m-2}_{i=0}\frac{f_{\delta^i}}{\mu^{1+i}}+\frac{A_2C_4}{\delta}\sum^{m-2}_{i=0}\frac{\sigma_{\delta^i}}{\mu^{1+i}}\bigg).
\end{equation*}

Let $\mu=\delta^\alpha~(\alpha>0)$, we calculate that
\begin{eqnarray*}
\sum_{i=1}^{m-2}\frac{f_{\delta^{i}}}{\mu^{i+1}}&=&\sum_{i=1}^{m-2}\frac{f_{\delta^{i}}}{\delta^{\alpha i+\alpha}}\\
&=&\frac{1}{1-\delta}\sum_{i=1}^{m-2}\frac{f_{\delta^{i}}}{\delta^{\alpha i+\alpha}\delta^{i-1}}(\delta^{i-1}-\delta^i)\\
&=&\frac{1}{(1-\delta)\delta^{2\alpha}}\sum_{i=1}^{m-2}\frac{f_{\delta^{i}}}{\delta^{(i-1)(1+\alpha)}}(\delta^{i-1}-\delta^i)\\
&\leq& \frac{1}{(1-\delta)\delta^{2\alpha}}\sum_{i=1}^{m-2}\int_{\delta^{i}}^{\delta^{i-1}}\frac{f_r}{r^{1+\alpha}}dr\\
&\leq& \frac{1}{(1-\delta)\delta^{2\alpha}}\int_{\delta^{m-1}}^{1}\frac{f_r}{r^{1+\alpha}}dr.
\end{eqnarray*}
Similarly,
\begin{equation*}
    \sum_{i=1}^{m-2}\frac{\sigma_{\delta^{i}}}{\mu^{i+1}}\leq \frac{1}{(1-\delta)\delta^{2\alpha}}\int_{\delta^{m-1}}^{1}\frac{\sigma_r}{r^{1+\alpha}}dr.
\end{equation*}
Therefore,
\begin{equation*}
K_{m}\leq \tilde C(\delta^m)^{\alpha}\bigg(1+\int_{\delta^{m-1}}^{1}\frac{f_r+\sigma_r}{r^{1+\alpha}}dr\bigg).
\end{equation*}
Then from Lemma \ref{ecl3.3}, we have that
\begin{equation*}
    u(x)\leq \tilde C\bigg\{(\delta^m)^{\alpha}\bigg(1+\int_{\delta^{m}}^{1}\frac{f_r+\sigma_r}{r^{1+\alpha}}dr\bigg)x_n+\delta^{m-1}(f_{\delta^{m-1}}+
    \sigma_{\delta^{m-1}})\bigg\}~~~~\mbox{in}~~\Omega_{\delta^m},
\end{equation*}
for $m=0,1,2,...$.

Suppose $r\leq \delta^2$. Let $m\geq 2$ such that $\delta^{m+1}<r\leq\delta^m$. Then for $x\in\Omega_r$, we have that
\begin{eqnarray*}
  u(x) &\leq& \tilde C\bigg\{(\delta^m)^{\alpha}\bigg(1+\int_{\delta^{m}}^{1}\frac{f_t+\sigma_t}{t^{1+\alpha}}dt\bigg)x_n+\delta^{m-1}(f_{\delta^{m-1}}+
    \sigma_{\delta^{m-1}})\bigg\}\\
   &\leq&  \tilde C\bigg\{\frac{r^{\alpha}}{\delta^{\alpha}}\bigg(1+\int_{r}^{1}\frac{f_t+\sigma_t}{t^{1+\alpha}}dt\bigg)x_n+
   \frac{r}{\delta^2}\big(f_{\frac{r}{\delta^2}}+\sigma_{\frac{r}{\delta^2}}\big)\bigg\}\\
   &\leq& Cr\bigg\{r^{\alpha}\bigg(1+\int_{r}^{1}\frac{f_t+\sigma_t}{t^{1+\alpha}}dt\bigg)+
   f_{\Lambda r}+\sigma_{\Lambda r}\bigg\},
\end{eqnarray*}
where $\Lambda=\frac{1}{\delta^2}$.\qed

\subsection{Differentiability at flat points}
In this subsection, we prove Theorem \ref{t1.1} with respect to flat points. We denote $\partial\Omega$ nearby 0 by
$$\partial\Omega\cap\{x||x'|\leq 1\}=\{(x',\nu(x'))||x'|\leq 1\},$$
where $\nu:R^{n-1}\rightarrow R$. From Lemma \ref{l2.3}, we know that $\nu$ is differentiable. Combining with $\eta=(0,...,0,1)$, we have that $\nabla \nu(0)=0$. Hence,
$$\frac{\nu(x')}{|x'|}\rightarrow 0 \mbox{ as } |x'|\rightarrow 0.$$
Let
$$\tilde\gamma(r)=\sup\{|\nu(x')|~|~|x'|\leq r\}$$
and
$$D(r)=\max\{\tilde\gamma(r),\gamma(r)\}.$$
Hence, $\frac{D(r)}{r}\rightarrow 0$ as $r\rightarrow 0$. Then just like (5) in Normalization \ref{nor1}, we assume that
$$D(t r)<tC,~\forall~ 0<t\leq 1,$$
where $C$ is the same constant in (5).

Under above assumptions, we observe that $\frac{\partial u(0)}{\partial x_i}=0$ for $i=1,2,...,n-1$ and $\frac{\partial u(0)}{\partial x_n}\geq 0.$ Hence, combining with Normalization \ref{nor1}, Theorem \ref{t1.1} is equivalent to:
\begin{theorem}\label{tt3.1}
There exists $a\geq 0$ such that
\begin{equation}\label{f}
|u(x)-ax_n|=o(|x|) ~as ~|x|\rightarrow 0.
%\frac{|u(x)-ax_n|}{|x|}\rightarrow 0 ~as ~|x|\rightarrow 0.
\end{equation}
\end{theorem}

\begin{remark}
Notice here that we only know that $\frac{|u(x)-ax_n|}{|x|}\rightarrow 0$ as $|x|\rightarrow 0$ and have no estimate of the convergence rate.
We prove (\ref{f}) by bounding the graph of $u$ between two hyperplanes in each scale and then use the principle that any bounded sequence has an accumulation to show the convergence of the slopes of the hyperplanes. Therefore no convergence rate can be obtained. Actually, Counterexample 4.1 and 4.3 in \cite{LW} imply that there should be no convergence rate, otherwise the gradient of $u$ will be continuous. \qed
\end{remark}

\begin{lemma}[\textbf{Key iteration}]\label{l3.2}
There exist positive constants $\delta(<1)$, $\mu(<1)$, $M$, $A_{1}$ and $A
_{2}$ depending only on $n$ and $\lambda$ such that if
\begin{equation}\label{e3.3}
kx_{n}-b\leq u(x)\leq Kx_n+B~~~~in~~\Omega_1
\end{equation}
for some nonnegative constants $b$, $B$, $k$ and $K$, then there exist nonnegative constants $\widetilde{k}$
and $\widetilde{K}$ such that
\begin{eqnarray}\label{14}
\begin{split}
&\widetilde{k}x_n-A_{1}\|f\|_{L^{n}(\Omega_{1})}-A_{2}(K+k+b)D(1)\\
&\leq u(x)\leq\widetilde{K}x_n+A_{1}\|f\|_{L^{n}(\Omega_{1})}+A_{2}(K+k+B)\gamma(1)~~~~in~~\Omega_{\delta},
\end{split}
\end{eqnarray}
where either
\begin{equation}\label{e3.4}
\widetilde{k}=\big(k-Mb+\mu(K-k)\big)^+~~\mbox{and}~~\widetilde{K}=K+MB,
\end{equation}
or\\
\begin{equation}\label{e3.5}
\widetilde{k}=(k-Mb)^+~~\mbox{and}~~\widetilde{K}=K+MB-\mu(K-k).
\end{equation}
\end{lemma}

\begin{remark}
From the geometric point of view, (\ref{e3.3}) means that the graph of $u$ lies between two hyperplanes, and (\ref{e3.4}) and (\ref{e3.5}) imply $\widetilde{K}-\widetilde{k}\leq (1-\mu)({K}-{k})+M(B+b)$, which means that the two hyperplanes approach to each other as the scale decreases.\qed
\end{remark}

\noindent\textbf{Proof.} 1. \emph{Claim:} There exist positive constants $M$, $\delta_1$ and $C_1$ depending only on $n$ and $\lambda$
such that
\begin{eqnarray}\label{17}
\begin{split}
&(k-Mb)x_{n}-C_1||f||_{L^{n}(\Omega_1)}-(k+Mb)D(1)\\
&\leq u(x)\leq (K+MB)x_n+C_1\|f\|_{L^{n}(\Omega_{1})}+(K+MB)\gamma(1)~~~~\mbox{in}~~\Omega_{\delta_1}
\end{split}
\end{eqnarray}
\proof Let $M_1=2\sqrt{n-1}\left(1+\frac{2\sqrt{n-1}}{\lambda}\right)$, $M=4M_1$, $\delta_1=\frac{1}{2M_1}$ and $\Psi=\Psi_{\frac{1}{2},\delta_1}$ (see Lemma \ref{l3.1} and recall $\gamma(1)\leq D(1)\leq\frac{\delta_1}{4}$). Let $v(x)=u(x)-Kx_n-B\Psi(x)$. We claim that
\begin{equation*}
    v(x)\leq K\gamma(1)~~~~~~~\mbox{on}~~\partial\Omega[\frac{1}{2},\delta_1].
\end{equation*}
In fact, we separate $\partial\Omega[\frac{1}{2},\delta_1]$ into two parts: $\partial Q[\frac{1}{2},\delta_1]\cap\bar\Omega$ and $\partial\Omega\cap \overline{ Q[\frac{1}{2},\delta_1]}$. On the first part, since $\Psi\geq 1$, we have $v\leq 0$. On the second part, since $\Psi\geq 0$, $u=0$ and $x_n\geq -\gamma(1)$, we have $v\leq K\gamma(1)$. Thus, we have proved the claim.
Then it is easy to verify that
\begin{eqnarray*}
 \left\{\begin{array}{l}
 \displaystyle-a^{ij}(x)\frac{\partial^2v(x)}{\partial x_i \partial x_j}\leq f(x)~~~~\mbox{in}~~\Omega[\frac{1}{2},\delta_1];\\
v\leq K\gamma(1)~~~~\mbox{on}~~\partial\Omega[\frac{1}{2},\delta_1].\\
\end{array}
\right.
\end{eqnarray*}

According to the Aleksandrov-Bakelman-Pucci maximum principle\cite{GT}, there exists a positive constant $C_1$
depending only on $n$ and $\lambda$ such that
\begin{equation*}
v(x)\leq K\gamma(1)+C_1||f||_{L^{n}(\Omega_1)}~~~~\mbox{in}~~\Omega[\frac{1}{2},\delta_1].
\end{equation*}
From (\ref{11}(4)) and the definition of $v$, we have that
\begin{eqnarray*}\label{18}
u(x)&=&Kx_n+B\Psi(x)+v(x)\\
&\leq& (K+\frac{2B}{\delta_1})x_n+C_1||f||_{L^{n}(\Omega_1)}+(K+\frac{2B}{\delta_1})\gamma(1)~~~~\mbox{in}~~\Omega_{\delta_1},
\end{eqnarray*}
or (recall $M=\frac{2}{\delta_1}$)
\begin{equation*}\label{19}
u(x)\leq (K+MB)x_n+C_1||f||_{L^{n}(\Omega_1)}+(K+MB)\gamma(1)~~~~\mbox{in}~~\Omega_{\delta_1}.
\end{equation*}
Therefore, we have proved the right side of (\ref{17}).

On the other hand, let $w(x)=kx_n-b\Psi(x)-u(x)$. Similarly, We claim that
\begin{equation*}\label{20}
    w(x)\leq kD(1)~~~~~~~\mbox{on}~~\partial\Omega[\frac{1}{2},\delta_1].
\end{equation*}
In fact, we separate $\partial\Omega[\frac{1}{2},\delta_1]$ into two parts: $\partial Q[\frac{1}{2},\delta_1]\cap\bar\Omega$ and $\partial\Omega\cap \overline{Q[\frac{1}{2},\delta_1]}$. On the first part, since $\Psi\geq 1$, we have $w\leq 0$. On the second part, since $\Psi\geq 0$, $u=0$ and $x_n\leq D(1)$, we have $w\leq kD(1)$. Thus, we have proved the claim. Then it is easy to verify that (recall $f\geq 0$)
\begin{eqnarray*}
 \left\{\begin{array}{l}
\displaystyle-a^{ij}(x)\frac{\partial^2w(x)}{\partial x_i \partial x_j}\leq 0~~~~\mbox{in}~~\Omega[\frac{1}{2},\delta_1];\\
w\leq kD(1)~~~~\mbox{on}~~ \partial\Omega[\frac{1}{2},\delta_1].\\
\end{array}
\right.
\end{eqnarray*}

According to the Aleksandrov-Bakelman-Pucci maximum principle, we have that
\begin{equation*}\label{21}
w(x)\leq kD(1)~~~~\mbox{in}~~  \Omega[\frac{1}{2},\delta_1].
\end{equation*}
Similarly, from (\ref{11}(4)) and the definition of $w$, we have that
\begin{eqnarray*}
u(x)&=&kx_n-b\Psi(x)-w(x)\\
  &\geq& kx_n-\frac{2\big(x_n+\gamma(1)\big)}{\delta_1}b-kD(1) \\
  &=&(k-Mb)x_n-(k+Mb)D(1)~~~~\mbox{in}~~\Omega_{\delta_1}.
\end{eqnarray*}
Hence, the left side of (\ref{17}) holds. Therefore, we have proved Claim 1.

2. Let $\delta=\frac{\delta_1}{2M_1}$ and $\Gamma={\delta e_n+T_{M_1\delta}}$ (recall $D(1)<\frac{1}{4}\delta$). Then $\Gamma\subset\Omega_{\frac{1}{2}\delta_1}$ and $\mbox{dist}(\Gamma,\partial\Omega)\geq \frac{\delta}{2}$.
In the following, we will prove the lemma according to two cases: $u(\delta e_n)\geq \frac{1}{2}(K+k)\delta$ and $u(\delta e_n)\leq \frac{1}{2}(K+k)\delta$, corresponding to which (\ref{e3.4}) and (\ref{e3.5}) hold respectively.

\noindent\textbf{Case 1}: $u(\delta e_n)\geq \frac{1}{2}(K+k)\delta.$ Let $$v(x)=u(x)-(k-Mb)x_n+C_1||f||_{L^{n}(\Omega_1)}+(k+Mb)D(1).$$
Then
\begin{eqnarray*}
  v(\delta e_n)&=& u(\delta e_n)-(k-Mb)\delta+C_1||f||_{L^{n}(\Omega_1)}+(k+Mb)D(1)\\
&\geq& \left(\frac{1}{2}(K-k)+Mb\right)\delta+C_1||f||_{L^{n}(\Omega_1)}+(k+Mb)D(1).
\end{eqnarray*}
Since $v\geq0~~\mbox{in}~\Omega_{\delta_1}$, by the Harnack inequality\cite{GT}, we have that

\begin{equation*}\label{22}
    \sup_{\Gamma}v(x)\leq C_2\left(\inf_{\Gamma}v(x)+\|f\|_{L^{n}(\Omega_{1})}\right),
\end{equation*}
where $C_2(>3+C_1)$ is a constant depending only on $n,~\lambda,~\delta_1~\mbox{and}~\delta$. Therefore,
\begin{equation*}\label{23}
  \inf_{\Gamma}v(x)\geq\left\{\frac{1}{C_2}\left(\Big(\frac{1}{2}(K-k)+Mb\Big)\delta+(k+Mb)D(1)\right)+\Big(\frac{C_1}{C_2}-1\Big)\|f\|_{L^{n}(\Omega_{1})}\right\}^{+}:=a.
\end{equation*}

Let $\tilde{\Psi}=\tilde\Psi_{M_1\delta,\delta}$ and $w=a\tilde{\Psi}-v$. We claim that
\begin{equation*}\label{24}\
    w\leq 4M(K+k+b)D(1)~~~~~~\mbox{on}~~\partial\Omega[M_1\delta,\delta].
\end{equation*}
In fact, we separate $\partial\Omega[M_1\delta,\delta]$ into three parts:
$\{x\in\partial Q[M_1\delta,\delta]|x_n=\delta\}\cap\bar\Omega$, $\{x\in\partial Q[M_1\delta,\delta]|x_n\neq\delta\}\cap\bar\Omega$ and $\partial\Omega\cap \overline{Q[M_1\delta,\delta]}.$
On the first part, since $\tilde\Psi\leq 1$ and $v\geq a$, we have $w\leq 0$. On the second part, since $\tilde\Psi\leq 0$ and $v\geq 0$, we have $w\leq 0$. On the last part, since $v\geq 0~\mbox{and}~\tilde{\Psi}(x)\leq\frac{x_n+\gamma(2M_1\delta)}{\delta}$, we have that
\begin{eqnarray*}
    w(x)&=&a\tilde{\Psi}(x)-v(x)\leq  a\tilde{\Psi}(x)\\
    &\leq&\frac{1}{C_2}\Big((K+Mb)\delta+
    (k+Mb)D(1)\Big)\frac{x_n+\gamma(2M_1\delta)}{\delta}\\
 &\leq&\frac{1}{C_2}(K+Mb)\Big(x_n+\gamma(2M_1\delta)\Big)+\frac{1}{C_2}(k+Mb)D(1)\frac{x_n+\gamma(2M_1\delta)}{\delta}.
\end{eqnarray*}
Since $\gamma(2M_1\delta)<\delta~,~C_2>1~\mbox{and}~|x_n|\leq D(1)<\delta$, we have that
\begin{equation*}\label{25}
   w\leq 4M(K+k+b)D(1)~~~~~\mbox{on}~~\partial\Omega\cap \overline{Q[M_1\delta,\delta]}.
\end{equation*}
Thus, we have proved the claim. Then it is easy to verify that
\begin{eqnarray*}
 \left\{\begin{array}{l}
\displaystyle-a^{ij}(x)\frac{\partial^2w(x)}{\partial x_i \partial x_j}\leq 0~~~~\mbox{in}~~\Omega[M_1\delta,\delta];\\
w\leq 4M(K+k+b)D(1)~~~~\mbox{on}~~\partial\Omega[M_1\delta,\delta].\\
\end{array}
\right.
\end{eqnarray*}

According to the Aleksandrov-Bakelman-Pucci maximum principle, we have that
\begin{equation*}\label{26}
w\leq 4M(K+k+b)D(1)~~~~\mbox{in}~~\Omega[M_1\delta,\delta].
\end{equation*}
From (\ref{12}(4)) and the definition of $v$, we have that
\begin{eqnarray*}
u(x)&=&(k-Mb)x_n-C_1||f||_{L^{n}(\Omega_1)}-(k+Mb)D(1)+v(x) \\
&=&(k-Mb)x_n-C_1||f||_{L^{n}(\Omega_1)}-(k+Mb)D(1)+a\tilde{\Psi}(x)-w(x)\\
&\geq&(k-Mb)x_n-C_1||f||_{L^{n}(\Omega_1)}-(k+Mb)D(1)+\\
&&\frac{a}{4\delta}\Big(x_n+\gamma(2M_1\delta)\Big)-4M(K+k+b)D(1)\\
&\geq&  (k-Mb)x_n-C_1||f||_{L^{n}(\Omega_1)}-5M(K+k+b)D(1)+\\
&&\frac{x_n+\gamma(2M_1\delta)}{4\delta}\bigg(\frac{K-k}{2C_2}\delta-\|f\|_{L^{n}(\Omega_{1})}\bigg)\\
&\geq& \bigg(\frac{K-k}{8C_2}+k-Mb\bigg)x_n-(C_1+1)\|f\|_{L^{n}(\Omega_{1})}-\\
&&5M(K+k+b)D(1)+\frac{\gamma(2M_1\delta)}{8C_2}(K-k)\\
&\geq& \left(\frac{K-k}{8C_2}+k-Mb\right){x_n}-(C_1+1)\|f\|_{L^{n}(\Omega_{1})}-\\
&&(5M+1)(K+k+b)D(1)~~~~\mbox{in}~~\Omega_\delta.
\end{eqnarray*}
Let
\begin{equation}\label{e3.6}
\mu=\frac{1}{8C_2},~A_1=1+2C_1~\mbox{and}~A_2=6M+1.
\end{equation}
Then combining with (\ref{17}) and $u\geq 0$ we obtain that (\ref{14}) and (\ref{e3.4}) hold.

\noindent\textbf{Case 2}: $u(\delta e_n)\leq \frac{1}{2}(K+k)\delta.$ The proof is similar to Case 1. Let
$$v(x)=(K+MB)x_n+C_1\|f\|_{L^{n}(\Omega_{1})}+(K+MB)\gamma(1)-u(x).$$
Then
\begin{eqnarray*}
  v(\delta e_n)&=& (K+MB)\delta+C_1\|f\|_{L^{n}(\Omega_{1})}+(K+MB)\gamma(1)-u(\delta e_n)\\
&\geq& \left(\frac{1}{2}\left(K-k\right)+MB\right)\delta+C_1\|f\|_{L^{n}(\Omega_{1})}+(K+MB)\gamma(1).
\end{eqnarray*}
Since $v\geq0~~\mbox{in}~~ \Omega_{\delta_1}$, by the Harnack inequality, we have that

\begin{equation*}
    \sup_{\Gamma}v(x)\leq C_2\left(\inf_{\Gamma}v(x)+\|f\|_{L^{n}(\Omega_{1})}\right).
\end{equation*}
Hence,
\begin{equation*}\label{29}
  \inf_{\Gamma}v(x)\geq\left\{\frac{1}{C_2}\bigg(\Big(\frac{1}{2}(K-k)+MB\Big)\delta+(K+MB)\gamma(1)\bigg)+\Big(\frac{C_1}{C_2}-1\Big)\|f\|_{L^{n}(\Omega_{1})}\right\}^{+}:=a.
\end{equation*}

Let $\tilde{\Psi}=\tilde\Psi_{M_1\delta,\delta}~\mbox{and}~w=a\tilde{\Psi}-v$. We claim that
\begin{equation*}\label{30}\
    w\leq 0~~~~~~\mbox{on}~~\partial\Omega[M_1\delta,\delta].
\end{equation*}
In fact, we separate $\partial\Omega[M_1\delta,\delta]$ into three parts:
$\{x\in\partial Q[M_1\delta,\delta]|x_n=\delta\}\cap\bar\Omega$, $\{x\in\partial Q[M_1\delta,\delta]|x_n\neq\delta\}\cap\bar\Omega$ and $\partial\Omega\cap \overline{Q[M_1\delta,\delta]}.$ On the first part, since $\tilde\Psi\leq 1$ and $v\geq a$, we have $w\leq 0$. On the second part, since $\tilde\Psi\leq 0$ and $v\geq 0$, we have $w\leq 0$. On the last part, since $\tilde{\Psi}(x)\leq\frac{x_n+\gamma(2M_1\delta)}{\delta}$ and $C_2>(3+C_1)$, we have that
\begin{eqnarray*}
  a\tilde{\Psi}(x)-v(x)&\leq& \frac{1}{C_2}\Big((K+MB)\delta+
    (K+MB)\gamma(1)\Big)\frac{x_n+\gamma(2M_1\delta)}{\delta}-v(x)\\
 &\leq& \frac{1}{C_2}(K+MB)\Big(x_n+\gamma(2M_1\delta)\Big)+\frac{1}{C_2}(K+MB)\gamma(1)
 \frac{x_n+\gamma(2M_1\delta)}{\delta}\\
 &&-\bigg((K+MB)\Big(x_n+\gamma(1)\Big)+C_1\|f\|_{L^{n}(\Omega_{1})}-u(x)\bigg)\\
 &\leq& u(x)=0
  ~~~~\mbox{on}~~ \partial\Omega\cap \overline{Q[M_1\delta,\delta]}.
\end{eqnarray*}
Thus, we have proved the claim. Then it is easy to verify that
\begin{eqnarray*}
 \left\{\begin{array}{l}
\displaystyle-a^{ij}(x)\frac{\partial^2w(x)}{\partial x_i \partial x_j}\leq f(x)~~~~\mbox{in}~~\Omega[M_1\delta,\delta];\\
w\leq 0~~~~\mbox{on}~~ \partial\Omega[M_1\delta,\delta].\\
\end{array}
\right.
\end{eqnarray*}

According to the Aleksandrov-Bakelman-Pucci maximum principle, we have that

\begin{equation*}
w\leq C_1\|f\|_{L^{n}(\Omega_{1})}~~~~\mbox{in}~~ \Omega[M_1\delta,\delta].
\end{equation*}
From (\ref{12}(4)) and the definition of $v$, we have that
\begin{eqnarray*}
  u(x)&=&(K+MB)x_n+C_1||f||_{L^{n}(\Omega_1)}+(K+MB)\gamma(1)-v(x) \\
  &=&(K+MB)x_n+C_1||f||_{L^{n}(\Omega_1)}+(K+MB)\gamma(1)-a\tilde\Psi(x)+w(x)\\
&\leq& (K+MB)x_n+C_1||f||_{L^{n}(\Omega_1)}+(K+MB)\gamma(1)-\\
&&\frac{a}{4\delta}\Big(x_n+\gamma(2M_1\delta)\Big)
+C_1\|f\|_{L^{n}(\Omega_{1})}\\
&\leq& (K+MB)x_n+2C_1||f||_{L^{n}(\Omega_1)}+(K+MB)\gamma(1)-\\
&&\frac{x_n+\gamma(2M_1\delta)}{4\delta}\left(\frac{K-k}{2C_2}\delta-||f||_{L^{n}(\Omega_1)}\right)\\
&\leq& \Big(K+MB-\frac{ K-k}{8C_2}\Big)x_n+(2C_1+1)||f||_{L^{n}(\Omega_1)}+\\
&&(K+MB)\gamma(1)-\frac{\gamma(2M_1\delta)}{8C_2}(K-k)\\
&\leq& \Big(K+MB-\frac{ K-k}{8C_2}\Big)x_n+(2C_1+1)||f||_{L^{n}(\Omega_1)}+\\
&&( 2M+1)(K+k+B)\gamma(1)~~~~\mbox{in}~~\Omega_\delta.
\end{eqnarray*}
Combining with (\ref{17}), (\ref{e3.6}) and $u\geq 0$, we have that (\ref{14}) and (\ref{e3.5}) hold.\qed

As the proof of Lemma \ref{ecl3.3}, by $0\leq u\leq 1$ and Lemma \ref{l3.2}, we have
\begin{lemma}\label{l3.3}
There exist nonnegative sequences
$\{b_m\}_{m=0}^{\infty}$, $\{B_m\}_{m=0}^{\infty}$, $\{k_m\}_{m=0}^{\infty}$ and $\{K_m\}_{m=0}^{\infty}$ with $b_0=k_0=K_0=0$, $B_0=1$ and for
$m=0,1,2,...,$
\begin{eqnarray}\label{31}
 \left\{\begin{array}{l}
\displaystyle b_{m+1}=A_1\delta^m||f||_{L^{n}(\Omega_{\delta^m})}+A_2(K_m+k_m+\frac{b_m}{\delta^m})\delta^mD(\delta^m),\\[5pt]
\displaystyle B_{m+1}=A_1\delta^m||f||_{L^{n}(\Omega_{\delta^m})}+A_2(K_m+k_m+\frac{B_m}{\delta^m})\delta^m\sigma(\delta^m),\\[5pt]
 \mbox{and either}\\
 \displaystyle k_{m+1}=\big(k_m-M\frac{b_m}{\delta^m}+\mu(K_m-k_m)\big)^+\ \ \mbox{and}\ \ K_{m+1}=K_m+M\frac{B_m}{\delta^m},\\[5pt]
 \mbox{or}\\
 \displaystyle k_{m+1}=(k_m-M\frac{b_m}{\delta^m})^+\ \ \mbox{and}\ \ K_{m+1}=K_m+M\frac{B_m}{\delta^m}-\mu(K_m-k_m),
\end{array}
\right.
\end{eqnarray}
such that
\begin{equation}\label{e3.7}
k_mx_n-b_m\leq u(x)\leq K_mx_n+B_m\ \ in\ \ \Omega_{\delta^m}.
\end{equation}
\end{lemma}
~\\

\noindent\textbf{Proof of Theorem \ref{tt3.1}.} Let $\{b_m\}_{m=0}^{\infty}$,
$\{B_m\}_{m=0}^{\infty}$, $\{k_m\}_{m=0}^{\infty}$ and
$\{K_m\}_{m=0}^{\infty}$ be defined as in Lemma \ref{l3.3}. For
simplicity, we denote $||f||_{L^{n}(\Omega_r)}$ by $f_r$, $D(r)$ by $D_r$ and
$\sigma(r)$ by $\sigma_r$. We prove the theorem by the following several claims.

\noindent1. \emph{Claim 1:} $\{K_m+k_m\}_{m=0}^{\infty}$ is bounded.

\proof By induction, we have that
\begin{equation*}
    K_m\geq k_m \mbox{ for any }m\geq 0.
\end{equation*}
Since $K_{m+1}\leq K_m+M\frac{B_m}{\delta^m}$ (see (\ref{31})), we have that
\begin{eqnarray*}
   &&K_{m+1}+k_{m+1}+\frac{B_{m+1}}{\delta^{m+1}}\\
   &&\leq 2(K_m+M\frac{B_m}{\delta^m})+\frac{B_{m+1}}{\delta^{m+1}}\\
  &&\leq 2M\sum_{i=0}^{m}\frac{B_i}{\delta^i}+\frac{B_{m+1}}{\delta^{m+1}} \\
   &&=\frac{2M}{\delta}\sum_{i=1}^{m}\Big(A_1f_{\delta^{i-1}}+A_2\big(K_{i-1}+k_{i-1}+\frac{B_{i-1}}{\delta^{i-1}}\big)\sigma_{\delta^{i-1}}\Big)+2M+\\
   &&\frac{1}{\delta}\Big(A_1f_{\delta^{m}}+A_2\big(K_m+k_m+\frac{B_m}{\delta^m}\big)\sigma_{\delta^{m}}\Big)\\
   &&\leq   2M+\frac{1}{\delta}A_1f_{\delta^{m}}+\frac{2MA_1}{\delta}\sum_{i=1}^{m}f_{\delta^{i-1}}+\frac{2MA_2}{\delta}\sum_{i=1}^{m+1}\big(K_{i-1}+k_{i-1}+\frac{B_{i-1}}{\delta^{i-1}}\big)\sigma_{\delta^{i-1}}.
\end{eqnarray*}
From (\ref{e3.11}), it is to see that $\{K_m+k_m+\frac{B_m}{\delta^m}\}_{m=0}^{\infty}$ is bounded. Hence, $\{K_m+k_m\}_{m=0}^{\infty}$ is bounded.
~\\

\noindent2. \emph{Claim 2:} $\lim_{m\rightarrow\infty}\frac{b_m}{\delta^m}=0$ and $\lim_{m\rightarrow\infty}\frac{B_m}{\delta^m}=0$.

\proof Since the proofs of $\{b_m\}_{m=0}^{\infty}$ and $\{B_m\}_{m=0}^{\infty}$ are similar, we only prove the first identity. By Claim 1, there exists a positive constant $C_3$ such that
$0\leq K_m+k_m\leq C_3$ for any $m\geq 0$. Then combining with $f_{\delta^{m}}\leq 1$ and
$D(\delta^m)\leq \frac{1}{A_2}$ (for an integer $m$ large enough), we have that
$$b_{m+1}=A_1\delta^m f_{\delta^{m}}+A_2\big(K_m+k_m+\frac{b_m}{\delta^m}\big)\delta^mD(\delta^m)\leq
A_1\delta^m+A_2C_3\delta^m+b_m.$$
Then $\{b_m\}$ is bounded since
$0<\delta<1.$ Since $$\frac{b_{m+1}}{\delta^{m+1}}=\frac{1}{\delta}A_1f_{\delta^{m}}+\frac{1}{\delta}A_2\big(K_m+k_m+\frac{b_m}{\delta^m}\big)D(\delta^m),$$
and $\{b_m\}$ is bounded, we conclude that $\{\frac{b_{m}}{\delta^{m}}\}$ is bounded (recall $D(\delta^m)\leq \delta^m$). Therefore, it follows that
$$\lim_{m\rightarrow\infty}\frac{b_{m+1}}{\delta^{m+1}}=
\frac{1}{\delta}\lim_{m\rightarrow\infty}\Big(A_1f_{\delta^{m}}+A_2(K_m+k_m+\frac{b_m}{\delta^{m}})D(\delta^m)\Big)=0.$$
~\\

\noindent3. \emph{Claim 3:} $\lim_{m\rightarrow\infty}(K_m-k_m)=0.$

\proof From $K_m\geq k_m$ and (\ref{31}), we have that
\begin{equation}\label{e3.12}
    0\leq K_{m+1}-k_{m+1}\leq (1-\mu)(K_m-k_m)+M\frac{b_m+B_m}{\delta^m},
\end{equation}
or
$$|K_{m+1}-k_{m+1}|\leq (1-\mu)|K_m-k_m|+\frac{M(2A_1+A_2C_4)}{\delta}(f_{\delta^{m-1}}+\sigma_{\delta^{m-1}}+D_{\delta^{m-1}}),$$
where $C_4$ is a positive constant such that $2(K_m+k_m+\frac{B_m}{\delta^{m}}+\frac{b_m}{\delta^{m}})\leq C_4$ for any $m\geq0$. The existence of such $C_4$ follows from Claim 1 and Claim 2. Let $C_5=\frac{M(2A_1+A_2C_4)}{\delta}$. It follows that
$$|K_{m+1}-k_{m+1}|\leq (1-\mu)^m\left(M+C_5\sum_{i=0}^{m-1}\frac{f_{\delta^{i}}+\sigma_{\delta^{i}}+D_{\delta^i}}{(1-\mu)^{i+1}}\right).$$
Let $1-\mu=\delta^\alpha~(\alpha>0)$, we calculate that
\begin{eqnarray*}
\sum_{i=1}^{m-1}\frac{f_{\delta^{i}}}{(1-\mu)^{i+1}}&=&\sum_{i=1}^{m-1}\frac{f_{\delta^{i}}}{\delta^{\alpha i+\alpha}}\\
&=&\frac{1}{1-\delta}\sum_{i=1}^{m-1}\frac{f_{\delta^{i}}}{\delta^{\alpha i+\alpha}\delta^{i-1}}(\delta^{i-1}-\delta^i)\\
&=&\frac{1}{(1-\delta)\delta^{2\alpha}}\sum_{i=1}^{m-1}\frac{f_{\delta^{i}}}{\delta^{(i-1)(1+\alpha)}}(\delta^{i-1}-\delta^i)\\
&\leq& \frac{1}{(1-\delta)\delta^{2\alpha}}\sum_{i=1}^{m-1}\int_{\delta^{i}}^{\delta^{i-1}}\frac{f_r}{r^{1+\alpha}}dr\\
&\leq& \frac{1}{(1-\delta)\delta^{2\alpha}}\int_{\delta^{m}}^{1}\frac{f_r}{r^{1+\alpha}}dr.
\end{eqnarray*}
For $D_{\delta^i}$ and $\sigma_{\delta^{i}}$, we have the similar result. Hence,
\begin{equation}\label{e3.13}
  |K_{m+1}-k_{m+1}|\leq C_6\delta^{m\alpha}(1+\int_{\delta^{m}}^{1}\frac{f_r+D_r+\sigma_r}{r^{1+\alpha}}dr).
\end{equation}
By L'Hospital rule, the right-hand side of (\ref{e3.13}) tends to 0 as $m\rightarrow \infty$. Thus, we have proved the claim.
~\\

\noindent4. \emph{Claim 4:} $\{K_m+k_m\}_{m=0}^{\infty}$ is convergent and we
set $$\lim_{m\rightarrow\infty}\frac{K_m+k_m}{2}=a.$$

\proof For any $m_0\geq 2$, we define $\{\tilde b_m\}_{m=m_0}^{\infty}$, $\{\tilde B_m\}_{m=m_0}^{\infty}$, $\{\tilde k_m\}_{m=m_0}^{\infty}$ and $\{\tilde K_m\}_{m=m_0}^{\infty}$ as follows:
$$\tilde b_{m_0}=b_{m_0},~\tilde B_{m_0}=B_{m_0},~\tilde k_{m_0}=k_{m_0},~\tilde K_{m_0}=K_{m_0},$$
and
\begin{eqnarray*}
 \left\{\begin{array}{l}
\displaystyle \tilde b_{m+1}=A_1\delta^mf_{\delta^{m}},~\tilde B_{m+1}=A_1\delta^mf_{\delta^{m}}+C_7\delta^m\sigma_{\delta^{m}},\\[5pt]
\displaystyle \tilde k_{m+1}=\big(\tilde k_m-M\frac{\tilde b_m}{\delta^m}+\mu(\tilde K_m-\tilde k_m)\big)^+\ \ \mbox{and}\ \ \tilde K_{m+1}=\tilde K_m+M\frac{\tilde B_m}{\delta^m}\\[5pt]
\displaystyle~~\mbox{if}~k_{m+1}=\big(k_m-M\frac{b_m}{\delta^m}+\mu(K_m-k_m)\big)^+\ \ \mbox{and}\ \ K_{m+1}=K_m+M\frac{B_m}{\delta^m},\\[5pt]
\displaystyle\tilde k_{m+1}=\big(\tilde k_m-M\frac{\tilde b_m}{\delta^m}\big)^+\ \ \mbox{and}\ \ \tilde K_{m+1}=\tilde K_m+M\frac{\tilde B_m}{\delta^m}-\mu(\tilde K_m-\tilde k_m)\\[5pt]
\displaystyle~~\mbox{if}~k_{m+1}=\big(k_m-M\frac{b_m}{\delta^m}\big)^+\ \ \mbox{and}\ \ K_{m+1}=K_m+M\frac{B_m}{\delta^m}-\mu(K_m-k_m),\\[5pt]
\end{array}
\right.
\end{eqnarray*}
where $C_7$ is a constant such that $A_2(K_m+k_m+\frac{B_m}{\delta^m})<C_7$ for any $m\geq 0$. Then $\tilde b_m\leq b_m$, $\tilde B_m\geq B_m$ for any $m\geq m_0$. Thus, by induction, we have that $\tilde k_m\geq k_m$, $\tilde K_m\geq K_m$ and $\tilde k_m\leq \tilde K_m$ for any $m\geq m_0$.

Likewise (\ref{e3.12}), we have that
\begin{equation*}
    0\leq \tilde K_{m+1}-\tilde k_{m+1}\leq (1-\mu)(\tilde K_m-\tilde k_m)+M\frac{\tilde b_m+\tilde B_m}{\delta^m},
\end{equation*}
or
\begin{equation}\label{e3.14}
    |\tilde K_{m+1}-\tilde k_{m+1}|\leq (1-\mu)|\tilde K_m-\tilde k_m|+\frac{M(2A_1+C_7)}{\delta}(f_{\delta^{m-1}}+\sigma_{\delta^{m-1}}).
\end{equation}
Let $1-\mu=\delta^\alpha~(\alpha>0)$, by the same arguments to derive (\ref{e3.13}), we have that
\begin{equation*}\label{e3.15}
    |\tilde K_{m+1}-\tilde k_{m+1}|\leq C_8\left(\delta^{(m-m_0)\alpha}|\tilde K_{m_0}-\tilde k_{m_0} |+\delta^{m\alpha}\int_{\delta^{m}}^{\delta^{m_0-2}}\frac{f_r+\sigma_r}{r^{1+\alpha}}dr\right).
\end{equation*}
By the definition of $\tilde K_m~\mbox{and}~\tilde k_m$, we have that
\begin{equation*}
    \tilde K_m+\tilde k_m-\mu(\tilde K_m-\tilde k_m)+M\frac{\tilde B_m-\tilde b_m}{\delta^m}\leq \tilde K_{m+1}+\tilde k_{m+1}\leq \tilde K_m+\tilde k_m+\mu(\tilde K_m-\tilde k_m)+M\frac{\tilde B_m}{\delta^m}.
\end{equation*}
Hence,
\begin{eqnarray*}\label{52}
|\tilde K_{m+1}+\tilde k_{m+1}-\tilde K_m-\tilde k_m|&\leq&\mu(\tilde K_m-\tilde k_m)+M\frac{\tilde B_m+\tilde b_m}{\delta^m}\nonumber\\
&\leq& \mu(\tilde K_m-\tilde k_m)+\frac{M(2A_1+C_7)}{\delta}(f_{\delta^{m}}+\sigma_{\delta^{m}}).
\end{eqnarray*}
Therefore,
\begin{eqnarray}\label{44}
  &&\sum_{j=m_0}^{\infty}|\tilde K_{j+1}+\tilde k_{j+1}-\tilde K_j-\tilde k_j|
  \leq \mu|\tilde K_{m_0}-\tilde k_{m_0}|+\frac{2M(A_1+C_7)}{\delta}\sum^{\infty}_{j=m_0-1}(f_{\delta^j}+\sigma_{\delta^j})\nonumber\\
  &&+\mu C_8\sum^{\infty}_{j=m_0}\left(\delta^{(j-m_0)\alpha}|\tilde K_{m_0}-\tilde k_{m_0} |+
  \delta^{j\alpha}\int_{\delta^{j}}^{\delta^{m_0-2}}\frac{f_r+\sigma_r}{r^{1+\alpha}}dr\right).
\end{eqnarray}

Now we estimate the right-hand side of (\ref{44}). Let $F_r:=\int_{r}^{\delta^{m_0-2}}\frac{f_s}{s^{1+\alpha}}ds$ and then
\begin{eqnarray*}
  &&\sum^{\infty}_{j=m_0}\delta^{j\alpha}\int_{\delta^{j}}^{\delta^{m_0-2}}\frac{f_r}{r^{1+\alpha}}dr\\
  &&=\sum^{\infty}_{j=m_0}\delta^{j\alpha}F_{\delta^{j}}\frac{\delta^{j}-\delta^{j+1}}{\delta^{j}-\delta^{j+1}}\\
  &&\leq\frac{1}{(1-\delta)\delta^{\alpha}}\sum^{\infty}_{j=m_0}\int_{\delta^{j+1}}^{\delta^{j}}r^{\alpha-1}F_rdr\\
  &&=\frac{1}{(1-\delta)\delta^{\alpha}}\int_{0}^{\delta^{m_0}}r^{\alpha-1}F_rdr\\
  &&=\frac{1}{(1-\delta)\delta^{\alpha}}\int_{0}^{\delta^{m_0}}r^{\alpha-1}\int_{r}^{\delta^{m_0-2}}\frac{f_s}{s^{1+\alpha}}dsdr\\
  &&=\frac{1}{(1-\delta)\delta^{\alpha}}\left(\int_{0}^{\delta^{m_0}}\frac{f_s}{s^{1+\alpha}}\int_{0}^{s}r^{\alpha-1}drds
  +\int_{\delta^{m_0}}^{\delta^{m_0-2}}\frac{f_s}{s^{1+\alpha}}\int_{0}^{\delta^{m_0}}r^{\alpha-1}drds\right)\\
  &&=\frac{1}{\alpha(1-\delta)\delta^{\alpha}}\left(\int_{0}^{\delta^{m_0}}\frac{f_r}{r}dr
  +(\delta^{m_0})^{\alpha}\int_{\delta^{m_0}}^{\delta^{m_0-2}}\frac{f_r}{r^{1+\alpha}}dr\right).
\end{eqnarray*}
Similarly,
\begin{eqnarray*}
    \sum^{\infty}_{j=m_0}\delta^{j\alpha}\int_{\delta^{j}}^{\delta^{m_0-2}}\frac{\sigma_r}{r^{1+\alpha}}dr
    =\frac{1}{\alpha(1-\delta)\delta^{\alpha}}\left(\int_{0}^{\delta^{m_0}}\frac{\sigma_r}{r}dr
  +(\delta^{m_0})^{\alpha}\int_{\delta^{m_0}}^{\delta^{m_0-2}}\frac{\sigma_r}{r^{1+\alpha}}dr\right).
\end{eqnarray*}
On the other hand, by the same arguments to derive (\ref{e3.11}), we have that
\begin{equation*}
    \sum^{\infty}_{j=m_0-1}f_{\delta^j}\leq \frac{\delta}{1-\delta}\int_{0}^{\delta^{m_0-2}}\frac{f_r}{r}dr~~\mbox{and~} \sum^{\infty}_{j=m_0-1}\sigma_{\delta^j}\leq \frac{\delta}{1-\delta}\int_{0}^{\delta^{m_0-2}}\frac{\sigma_r}{r}dr
\end{equation*}

Therefore,
\begin{eqnarray}\label{45}
 &&\sum_{j=m_0}^{\infty}|\tilde K_{j+1}+\tilde k_{j+1}-\tilde K_j-\tilde k_j|\nonumber\\
 &&\leq C_9\left(|\tilde K_{m_0}-\tilde k_{m_0}|+\int_{0}^{\delta^{m_0-2}}\frac{f_r+\sigma_r}{r}dr\right)+\nonumber\\
 &&C_9(\delta^{m_0})^{\alpha}\int_{\delta^{m_0}}^{\delta^{m_0}-2}\frac{f_r+\sigma_r}{r^{1+\alpha}}dr,
\end{eqnarray}
where $C_9$ is a positive constant.

We suppose by the contradiction that $\{k_m+K_m\}^{\infty}_{m=0}$ is not convergent. Since $\{k_m+K_m\}^{\infty}_{m=0}$ is bounded,
it has at least two different accumulation points which we denote by $0\leq \beta_1<\beta_2$. Since $\lim_{m\rightarrow\infty}(K_m-k_m)=0$, $\tilde k_{m_0}=k_{m_0}$ and $\tilde K_{m_0}=K_{m_0}$, the right side of above inequality tends to zero as $m_0\rightarrow \infty.$ Then, there exists an integer $m_0$ large enough such that
\begin{equation*}\label{46}
    |K_{m_0}+k_{m_0}-\beta_1|<\frac{\beta_2-\beta_1}{4},
\end{equation*}
and
\begin{equation*}\label{47}
    \sum_{j=m_0}^{\infty}|\tilde K_{j+1}+\tilde k_{j+1}-\tilde K_j-\tilde k_j|<\frac{\beta_2-\beta_1}{4}.
\end{equation*}
Hence, for any $m\geq m_0$,
\begin{eqnarray*}
   \beta_2-(k_m+K_m)&\geq& \beta_2-(\tilde k_m+\tilde K_m)\\
  &\geq& (\beta_2-\beta_1)-|\beta_1-(\tilde k_{m_0}+\tilde K_{m_0})|-|(\tilde k_{m}+\tilde K_{m})-(\tilde k_{m_0}+\tilde K_{m_0})|\\
  &\geq& (\beta_2-\beta_1)-\frac{\beta_2-\beta_1}{2}\geq \frac{\beta_2-\beta_1}{2}.
\end{eqnarray*}
Therefore, we get a contradiction with that $\beta_2$ is an accumulation point of $\{k_m+K_m\}^{\infty}_{m=0}$.
~\\

\noindent5. \emph{Claim 5:} Let $a$ be given in Claim 4. Then for each $m=0,1,2,...$, there exists $C_m$ such that
$\lim_{m\rightarrow\infty}C_m=0$ and $|u(x)-ax_n|\leq C_m\delta^m$ for any $x\in\Omega_{\delta^m}$.

\proof For any $m\geq0$ and any $x\in\Omega_{\delta^m}$, we have that
\begin{equation}\label{e3.16}
|u(x)-a x_n|\leq|u(x)-\frac{K_m+k_m}{2}x_n|+|(\frac{K_m+k_m}{2}-a)x_n|.
\end{equation}
From (\ref{e3.7}), $$-\frac{K_m-k_m}{2}x_n-b_m\leq
u(x)-\frac{K_m+k_m}{2}x_n\leq \frac{K_m-k_m}{2}x_n+B_m.$$ Let
$\tilde{C}_m=\max\{\frac{K_m-k_m}{2}+\frac{b_m}{\delta^m},\frac{K_m-k_m}{2}+\frac{B_m}{\delta^m}\}$
and then
\begin{equation}\label{e3.17}
|u(x)-\frac{K_m+k_m}{2}x_n|\leq \tilde{C}_m\delta^m.
\end{equation}
From Claim 2 and Claim 3, we have that
\begin{equation}\label{e3.18}
\lim_{m\rightarrow\infty}\tilde{C}_m=0.
\end{equation}
Let $C_m=\tilde{C_m}+|\frac{K_m+k_m}{2}-a|$. From (\ref{e3.16}) and (\ref{e3.17}), it follows that
$$|u(x)-a x_n|\leq C_m\delta^m.$$
From Claim 4 and (\ref{e3.18}), we have that $C_m\rightarrow0$ as $m\rightarrow\infty$.

From Claim 5, we deduce that $u$ is differentiable at $0$ with derivative $ae_n$. Therefore, the proof of Theorem \ref{tt3.1} is completed.\qed

\begin{remark}\label{r3.1}
If  $u=g~\mbox{on}~\partial\Omega$ in (\ref{e1.1}) and $g$ satisfies the Dini condition (see (1.2) in \cite{LW2}),
it is not hard to verify that Theorem \ref{t1.1} is also true and the Dini condition can not be weaken. As pointed out in \cite{LW}, no continuity of the gradient of solutions can be expected. So far, together with related papers, we gain a deep insight into the boundary differentiability of solutions of (\ref{e1.1}).\qed
\end{remark}



%% References


\bibliographystyle{elsarticle-num}






% \end{thebibliography}

\end{document}